\documentclass[11pt]{article}

\usepackage{graphicx}

\usepackage{amssymb}

\begin{document}
\tolerance=10000
\tolerance=10000


\hyphenpenalty=2000
\hyphenation{visco-elastic visco-elasticity}
\setcounter{page}{1}
\thispagestyle{empty}




\font\note=cmr10 at 10 truept  


\def\pni{\par \noindent}
\def\vsh{\smallskip}
\def\s{\smallskip}
\def\vs{\medskip}
\def\vvs{\bigskip}
\def\vvvs{\bigskip\medskip} 
\def\vsp{\vsh \par\noindent}
\def\vsn{\vsh\pni}
\def\cen{\centerline}
\def\ra{\item{a)\ }} \def\rb{\item{b)\ }}   \def\rc{\item{c)\ }}
\def\eg{{\it e.g.}\ } \def\ie{{\it i.e.}\ }


\def\sg{\hbox{sign}\,}
\def\sgn{\hbox{sign}\,}
\def\sign{\hbox{sign}\,}
\def\e{{\rm e}}
\def\exp{{\rm exp}}
\def\ds{\displaystyle}
\def\dis{\displaystyle}
\def\q{\quad}    \def\qq{\qquad}
\def\lan{\langle}\def\ran{\rangle}
\def\l{\left} \def\r{\right}
\def\lra{\Longleftrightarrow}
\def\d{\partial}
\def\dr{\partial r}  \def\dt{\partial t}
\def\dx{\partial x}   \def\dy{\partial y}  \def\dz{\partial z}
\def\rec#1{{1\over{#1}}}
\def\zr{z^{-1}}



\def\hatt{\widehat}
\def\epsilons{{\widetilde \epsilon(s)}}
\def\sigmas{{\widetilde \sigma (s)}}
\def\fs{{\widetilde f(s)}}
\def\Js{{\widetilde J(s)}}
\def\Gs{{\widetilde G(s)}}
\def\Fs{{\wiidetilde F(s)}}
 \def\Ls{{\widetilde L(s)}}
\def\L{{\cal L}} 
\def\F{{\cal F}} 


\def\NN{{\bf N}}
\def\RR{{\bf R}}
\def\CC{{\bf C}}
\def\ZZ{{\bf Z}} 


\def\I{{\cal I}}  
\def\D{{\cal D}}  


\def\Gc{{\cal {G}}_c}   \def\Gcs{\barr{\Gc}} 
\def\Gs{{\cal {G}}_s}   \def\Gss{\barr{\Gs}} 
\def\Gck{\hatt{\Gc}} 
\def\args{(x/ \sqrt{D})\, s^{1/2}}
\def\argsa{(x/ \sqrt{D})\, s^{\beta}}
\def\arg{ x^2/ (4\,D\, t)}


\def\erf{\hbox{erf}}  \def\erfc{\hbox{erfc}}


\def\uks{{\widehat{\widetilde {u}}} (\kappa,s)}

\def\psikappa{\psi_\alpha^\theta(\kappa)}





\baselineskip=14pt  
\cen{FRACALMO PRE-PRINT   {\bf www.fracalmo.org}}
\vsh
\cen{To be published in}
\cen{\bf  Journal of Computational and Applied Mathematics}
\vsh
\hrule
   \vskip 1.0truecm

\hyphenpenalty=2000




\font\title=cmbx12 scaled\magstep2
\font\bfs=cmbx12 scaled\magstep1

\begin{center}

{\title Some recent advances}
\vs  

{\title in theory and simulation of}
\vs   

{\title fractional diffusion processes}

\vskip 0.5truecm

{Rudolf GORENFLO}$^{(1)}$
 and {Francesco MAINARDI} $^{(2)}$


\vs

$\null^{(1)}$
 Fachbereich Mathematik \& Informatik \\
Erstes Mathematisches Institut,
 Freie   Universit\"at  Berlin, \\
 Arnimallee  3, D-14195 Berlin, Germany \\
E-mail: {\tt gorenflo@mi.fu-berlin.de}
\\ [0.25 truecm]

$\null^{(2)}$
 Dipartimento di Fisica, Universit\`a di Bologna, and INFN, \\
Via Irnerio 46, I-40126 Bologna, Italy \\
E-mail: {\tt mainardi@bo.infn.it}
\\ [0.25 truecm]

{Revised Version : May 2008}
\end{center}

\noindent
{\small This paper is based on the invited lecture by Rudolf  Gorenflo 
at  the Second International Workshop on Analysis and Numerical Approximation 
of Singular Problems (IWANASP 2006) held  on 6-8 September 2006,
at the Aegean University in Karlovassi, Samos, Greece, 
jointly with the University of Portsmouth, UK and the Instituto Superior Tecnico (CEMAT), Lisbon, Portugal.}

\begin{abstract}
\noindent 
To offer an insight into the rapidly developing theory of fractional diffusion processes 
we describe in some detail three topics of current interest: 
(i) the well-scaled passage to the limit from continuous time random walk 
under power law assumptions to space-time fractional diffusion, 
(ii) the asymptotic universality of the Mittag-Leffler waiting time law in time-fractional processes, 
(iii)  our method of parametric subordination for generating particle trajectories.
\end{abstract}
\newpage
\def\sg{\hbox{sign}\,}
\def\sgn{\hbox{sign}\,}
\def\sign{\hbox{sign}\,}
\def\e{{\rm e}}
\def\exp{{\rm exp}}
\def\ds{\displaystyle}
\def\dis{\displaystyle}
\def\lan{\langle}\def\ran{\rangle}
\def\lt{\left} \def\rt{\right}  
\def\lra{\Longleftrightarrow}
\def\d{\partial}
\def\dr{\partial r}  \def\dt{\partial t}
\def\dx{\partial x}   \def\dy{\partial y}  \def\dz{\partial z}
\def\rec#1{{1\over{#1}}}
\def\zr{z^{-1}}
\def\hatt{\widehat}
\def\epsilons{{\widetilde \epsilon(s)}}
\def\sigmas{{\widetilde \sigma (s)}}
\def\fs{{\widetilde f(s)}}
\def\Js{{\widetilde J(s)}}
\def\Gs{{\widetilde G(s)}}
\def\Fs{{\wiidetilde F(s)}}
 \def\Ls{{\widetilde L(s)}}
\def\L{{\mathcal L}} 
\def\F{{\mathcal F}} 
\def\NN{{\bf N}}
\def\RR{{\bf R}}
\def\CC{{\bf C}}
\def\ZZ{{\bf Z}} 
\def\I{{\cal I}}  
\def\D{{\cal D}}  
\def\erf{\hbox{erf}}  \def\erfc{\hbox{erfc}}
\def\uks{{\widehat{\widetilde {u}}} (\kappa,s)}
\def\psikappa{\psi_\alpha^\theta(\kappa)}

\section{Introduction}  

The field of fractional (more generally anomalous) diffusion processes 
in recent decades has won more and more interest in applications in the sciences, in physics and chemistry, 
and even in finance. 
Here we will give  some insight into this rapidly developing field. 
Instead of trying to cover the whole spectrum of this field we  
will focus on three contrasting aspects to which we ourselves have affinity by our research, 
thereby trying  to keep this paper as self-contained as  possible. 
We hope not only to whet the appetite among people not fully
 familiar with the subject but also to
propagate our methodological viewpoints in the fractional diffusion
community. 
\vsp
Viewing fractional diffusion processes as generalization of the familiar classical diffusion process, 
we begin by considering the Cauchy problem for the classical diffusion equation
$$ \frac{\d u(x,t)}{\d t} = \frac{\d^2 u(x,t)}{\d x^2}\,, \q u(x,0^+)=f(x)\,, \q x\in \RR\,,\; t\ge 0\,.
\eqno(1.1)$$
In fractional diffusion equations (see Section 2 for explanations) the differentiations with respect 
to $t$  and $x$ are replaced by differentiations of non-integer order. 
For equation (1.1) it is well known that its solution $u(x,t)$ has the essential properties 
we expect from a diffusion process, 
that is a process of re-distribution in space $x$ and time $t$. 
Considering $u(x,t)$ as the spatial density of an extensive quantity, e.g.
mass, charge, or probability, we have {\bf (a)}, {\bf (b)} and {\bf (c)}:
\\ 
{\bf (a)} conservation of the total quantity:
$  \int_{-\infty}^{+\infty} \!\! u(x,t)\, dx = \int_{-\infty}^{+\infty}\!\! u(x,0^+)\, dx$, $\forall \,t>0$.
\\ 
{\bf (b)} preservation of non-negativity:
$ u(x,0^+) \ge 0$, $\forall \, x\in \RR$ implies 
   $u(x,t) \ge 0$, $\forall \,x\in \RR$, $t>0$.
\\
 {\bf (c)}  Another essential characteristic of problem (1.1) concerns
   the law of spreading (or dispersion) of the quantity. 
 With the special initial condition $u(x,0^+)= \delta(x)$
   (the Dirac generalized function), the variance grows linearly in time,
   that is 
$ \sigma^2(t) :=  \int_{-\infty}^{+\infty} x^2\, u(x,t)\, dx = 2t$.
More generally, in a classical diffusion process the variance,
which is a natural and common quadratic  measure of the spread
of a diffusing substance, grows {\it linearly} in time, that is,  
 if we allow a drift,  
we have  
$ \sigma^2(t) :=  \int_{-\infty}^{+\infty} x^2\, \left[u(x,t)-m(t)\right]\, dx 
 \sim   C\, t$ as $t \to\infty$
with
 $m(t):=  \int_{-\infty}^{+\infty} x\, u(x,t)\, dx$,
  for some constant $C>0$.
 \vsp
 The above properties {\bf (a)}, {\bf (b)} and {\bf (c)}
  are indeed shared by many processes governed by second-order 
linear parabolic equations.
Usurping the term {\it diffusion} for processes having properties {\bf (a)} and {\bf (b)}
but not necessarily {\bf (c)},
we follow the custom to call processes of {\it anomalous diffusion}
those in which, for initial condition $u(x,0^+)= \delta(x)$,
the variance does not exhibit essentially linear grow with $t\to \infty$.
Among these processes we single out the {\it sub-diffusive} ones for
    which the variance grows (for large $t$) more slowly than
    linearly,and the {\it super-diffusive} for which it grows (for large 
    $t$) faster than linearly, or even does not exist (i.e is infinite).
Focusing our attention to 
the space-time fractional diffusion equation (i.e. a {\it generalization} of the
classical diffusion equation (1.1) via suitable {\it pseudo-differential operators} 
interpreted as time and space {\it derivatives of fractional order}),
we will discuss the construction and  properties of its fundamental solution
(obtained for $f(x)= \delta(x)$) and its approximation by continuous time random walk.
This generalized diffusion equation has found wide interest among researchers in recent 20 years.
After a general survey of basic facts we will go into details
of three distinct but related topics.  
As this paper cannot be a substitute for an extensive monograph,
our presentation will naturally be biased by our and our close
co-workers' contributions.
We will meet the two complementary aspects of diffusive processes.
The {\it first} is the {\it macroscopic} aspect: the structure of the 
fundamental solutions in dependence on space $x$ and time $t$,
in particular their scaling properties and asymptotics.
Here $u(x,t)$  is viewed as the density with respect to $x$ evolving in $t$.
The {\it second} is the {\it microscopic} aspect:
the structure of the trajectories (paths) of particles subject to such process.
Here  $ u(x,t)$  is viewed as the sojourn probability density with respect to $x$ evolving in $t$. 
\vsp
The rough structure of our paper is as follows. In Section 2 we will give a survey
of the space-time fractional diffusion equation and
the essential properties of its fundamental solution.
Sections 3 and 4 are devoted to topic (i): {\it continuous  time random walk} (CTRW), 
Section 5 to topic (ii): {\it power laws} and 
{\it well-scaled passage to the diffusion limit}, 
Section 6 to topic (iii): our method of {\it parametric subordination} for 
exact simulation of  {trajectories}. 
Finally, in Section 7, we will draw some conclusions.  
\vsp
Of course, there are more significant advances than we can report here.
Let us only mention processes with distributed orders of fractional differentiation and
multi-dimensional processes and cite, e.g. the papers 
\cite{Andries-Umarov-Steinberg FCAA06,ChechkinGorenfloSokolov PRE02,ChechkinGorenfloSokolovGonchar FCAA03,
SokolovChechkinKlafter 04}   
and the fundamental monograph by Meerschaert \& Scheffler \cite{MMM-Scheffler BOOK01}. 
Much work has been done in numerical analysis of difference schemes for calculating
densities of fractional diffusion processes, see e.g. \cite{Ilic-et-al FCAA06}.
We apologize to all authors who feel that we appreciate their work insufficiently.  
Throughout this paper we will make liberal use of generalized functions
in the sense of Gel'fand and Shilov \cite{GelfandShilov 64}, 
so avoiding the cumbersome notations
of measures and functionals, and usually we will assume the occurring
(generalized) functions so well-behaved that our manipulations are allowed. 
By this we hope our presentation
to be welcome to researchers in applications as well as
 inspiring for pure mathematicians looking around
outside the ivory tower.
 \section{The space-time fractional diffusion}
 We begin by considering the Cauchy problem for the
(spatially one-dimensional) {\it space-time fractional diffusion equation}
$$  {\, _t}D_{*}^{\, \beta }\, u(x,t)
 \, = \,
 {\, _x}D_{ \theta}^{\,\alpha} \,u(x,t)\,,
\quad  u(x,0) = \delta (x)\,, \quad x \in \RR,\quad t \ge 0\,, \eqno(2.1) $$
where 
   $\alpha ,\,\theta,\, \beta $ are real parameters
 restricted to the ranges
$$ 0<\alpha\le 2\,,\quad  |\theta| \le \min \{\alpha, 2-\alpha\}\,,
  \quad 0<\beta\leq 1\,.\eqno(2.2) $$
Here
${ \,_t}D_*^{\,\beta}  $
denotes  the
{\it  Caputo-Dzherbashyan  fractional derivative}
of order $\beta $, acting on a sufficiently well-behaved function $f(t)$ of the time variable $t$,
$$
    _tD_*^\beta \,f(t) := 
       {\ds \rec{\Gamma(1-\beta )}}\,
	   {\ds\int_0^t
 {\ds {f^{(1)}(t') \over (t-t')^\beta}\, dt'
 }} \,,
 \q 0<\beta	<1\,, \eqno(2.3a) $$
and  $ {\,_x}D_{\,\theta}^{\,\alpha}$
denotes
the {\it  Riesz-Feller fractional derivative}
of order $\alpha $ and
skewness $\theta$,   acting on a sufficiently well-behaved function $g(x)$ of 
the space variable $x$,
$$\q\q  \,_xD_\theta^\alpha \,g(x)
 = {\Gamma(1+\alpha) \over \pi } \,
 \left\{\sin \,\left[(\alpha+\theta) {\pi\over 2}\right] \,
 \int_0^\infty
 {g(x+x')- g(x)  \over {x'}^{1+\alpha}}\, d x' \right.
\q\q\q\q \eqno(2.4a)$$
$$ + \left.\sin \,\left[(\alpha-\theta) {\pi \over 2}\right] \,
 \int_0^\infty
 {g(x-x')- g(x)  \over {x'}^{1+\alpha}}\, d x' \right\}\,, 
 \; 0<\alpha <2\,, \;  |\theta| \le \min \{\alpha, 2-\alpha\}\,.
 $$
In the symmetric case $\theta = 0$, which later will be 
     of  our main interest, formula (2.4a) simplifies to 
$$ \!\! \!_xD_0^\alpha \,f(x)
\! = \!{\Gamma(1+\alpha) \over \pi }
 \sin \left({\alpha \pi \over 2}\right) 
 \int_0^\infty 
 {g(x+x')- 2g(x) + g(x-x') \over {x'}^{1+\alpha}}\, d x'.
 \eqno(2.4'a)$$
The above representations of the space fractional derivatives are based on a suitable 
regularization of  hyper-singular integrals.
\vsp
In the limits $\beta = 1$  and $\alpha = 2$ (so $\theta=0$)
 we recover the first time derivative ${\ds {d f(t) \over dt}}$
and the second space derivative ${\ds {d^2 g(x) \over dx^2}}$, respectively.
\vsp
These representations mirror the fact that time-fractional (for $0<\beta<1$) 
processes are processes with {\it long memory} (see for this also Section 5)  
whereas space fractional (for $0<\alpha<2$) are processes with 
{\it spatial long-range interactions}.
For more information on these operators we refer to
\cite{GorMai CISM97,GorMai FCAA98,KST 06,Mainardi FCAA01,Podlubny 99,SKM 93}.
\vsp
Let us note that the solution $u(x,t)$  of the
 Cauchy problem (2.1), known as its {\it Green function} or fundamental solution,
can be viewed as the density of an extensive quantity or as
the probability density in the spatial variable $x$, evolving in time $t$.
In the case $\alpha =2$ (hence $\theta=0$) and $\beta =1$  we recover
the standard diffusion equation for which the fundamental solution
is the Gaussian density with variance $\sigma ^2 =2t$.
\vsp
Writing, with $\hbox{Re} [s] > \sigma _0$, $\kappa \in \RR$, the transforms of Laplace and Fourier as
$$     {\L} \lt\{ f(t);s\rt\}=  \widetilde f(s)
 := \int_0^{\infty} \e^{\ds \, -st}\, f(t)\, dt\,, $$
$$ {\F} \lt\{g(x);\kappa\rt\}=  \widehat g(\kappa)
  := \int_{-\infty}^{+\infty} \!\! \e^{\,\ds i\kappa x}\,g(x)\, dx\,,
$$
 the corresponding transforms
of ${ \,_t}D_*^{\,\beta} f(t)  $
and       $ {\,_x}D_{\,\theta}^{\,\alpha} g(x)$ are
$$   {\L} \lt\{  {\,_t}D_*^{\,\beta}\, f(t)\rt\} =
    s^{\, \ds \beta}\,  \widetilde{f} (s)- s^{\, \ds \beta -1}\, f(0)\,,
\eqno(2.3b)$$
$$\!\!   {\F} \lt\{ {\,_x}D_{\,\theta}^{\,\alpha}\, g(x )\rt\} =
  - |\kappa |^{\,\ds \alpha}  \, i ^{\,\ds \theta \,\sgn \kappa} \, \widehat {g}(\kappa)
  = - |\kappa |^{\,\ds \alpha}\,\e^ {\,\ds i \,(\sgn \kappa)\,\theta\,\pi/2 }\,  \widehat {g}(\kappa)\,.
   \eqno (2.4b)$$
\vsp
We will freely use the convolution theorems pertinent to these transforms, 
defining for generic functions: 
$(f_1 * f_2)(t)= \int_{[0,t]} f_1(t-t')\, f_2(t')\, dt', \;     t \ge 0$,
$ (g_1 * g_2)(x)= \int_{(-\infty,+\infty)} g_1(x-x')\, g_2(x')\, dx', \;   x \in \RR$,
and the convolution powers $f^{*n}(t)$  and  $g^{*n}(x)$ as $n$-fold convolutions ( $n \ge 0$).
 Note that $f^{*0} (t) = \delta(t)$, $g^{*0} (x) = \delta (x)$.  
For   mathematical details  
we cite
\cite{GorMai CISM97,Podlubny 99} on the Caputo-Dzherbashyan  derivative
and    \cite{SKM 93} on the Feller potentials.
For the general theory of pseudo-differential operators
and related Markov processes the interested reader is referred
to the excellent volumes by Jacob \cite{Jacob BOOKS}.
\vsp
Let us here recall the
representation in the Fourier-Laplace  domain of the (fundamental) solution
of (2.1).  
Using $\widehat \delta (\kappa ) \equiv 1$
we have from (2.1)
 $$  s^{\, \ds \beta}\,\widehat{\widetilde{u}}(\kappa ,s) - s^{\, \ds \beta -1}
    = -|\kappa|^{\, \ds \alpha} \,  i ^{\,\ds \theta \,\sgn \kappa}
   \, \widehat{\widetilde{u}}(\kappa ,s) \,,\eqno (2.5)
$$
hence explicitly 
$$  \widehat{\widetilde{u}}(\kappa ,s)
    =  \frac{ s^{\, \ds \beta -1}}
{s^{\,\ds \beta} + |\kappa |^{\, \ds \alpha} \,
i^{\,\ds \theta \,\sgn \kappa} }\,.
   \eqno(2.6)$$
For explicit expressions and plots of  the fundamental solution of (2.1)
in the space-time domain
we refer the reader to \cite{Mainardi FCAA01}.
There, starting from the fact  that the Fourier transform
$\widehat{u}(\kappa ,t)$ can be written as a Mittag-Leffler function
with complex argument, the authors
have  derived a Mellin-Barnes integral representation
of $u(x,t)$  with which they have proved the non-negativity
of the solution for values of the parameters
$\{\alpha,\, \theta, \,  \beta \}$ in the range (2.2)
and analyzed the evolution in time of its  moments.
In particular for $\{0<\alpha <2, \, \beta=1\}$ we obtain
the  densities of the stable processes of order $\alpha$ and skewness $\theta$.
The representation of $u(x,t)$ in terms of Fox $H$-functions
can be found  in   \cite{Mainardi JCAM05}.
We note, however,  that the solution of the space-time fractional diffusion equation (2.1) and
its variants has been investigated by several authors
as pointed out in the bibliography in \cite{Mainardi FCAA01}:
here we refer to some of them,
 \cite{BaeumerMeerschaert 01,Barkai PRE01,M3 PRE02sub,M3 PRE02sol,
Metzler-Klafter PhysRep00}, where the connection with the CTRW
was also pointed out.
\vsp
Henceforth our attention, if not said explicitly otherwise,
     will be focussed on the symmetric case $\theta = 0$. In this case
$u(x,t)$
is an even function of $x$ and  we get from (2.6)  the Fourier-Laplace representation 
$$  \widehat{\widetilde{u}}(\kappa ,s)
    =  \frac{ s^{\, \ds \beta -1}}
{s^{\,\ds \beta} + |\kappa |^{\, \ds \alpha}} \,,
   \eqno(2.7)$$
which allows  us demonstration of the conservation property (a), namely
$ \int_{-\infty}^{+\infty} \! u(x,t) \, dx \equiv 1$ for all $t>0$, 
and  calculation  
of the  variance  $\sigma^2(t) = \langle x^2(t)\rangle $ 
(the quadratic measure of the spread as  function of $t$)
 of a diffusing particle.
 \vsp
 From (2.7), more generally already from (2.6), 
 we find by aid of well-known properties of the Fourier transform
 $\widehat{\widetilde{u}}(0,s)= 1/s$, hence 
 ${\int_{-\infty}^{+\infty} u(x,t) \, dx} = \widehat{u}(0,t)= 1$
     and so  
	 ${ \int_{-\infty}^{+\infty} u(x,t) \, dx} \equiv 1$.
\vsp	 
 For the variance
  we find
 $$\sigma^2(t)=\int_{-\infty}^{+\infty}x^{2}\,u(x,t)\,dx 
     =
  - \left.\frac{\partial^{2}}{\partial\kappa^{2}}\widehat{u}(\kappa,t)\right|_{\kappa=0}\,.
  \eqno(2.8)$$
   Using the
  Mittag-Leffler function  
$$E_{\beta}(z)=\sum_{n=0}^{\infty} \frac{z^{n}}{\Gamma(1+n \beta)}\,,\eqno(2.9)$$  
 see \eg \cite{GorMai CISM97,Podlubny 99},
  we find by Laplace inversion the convergent series
$$\widehat{u}(\kappa,t)= E_{\beta}(-|\kappa|^{\alpha}\, t^{\beta})=1-
\frac{|\kappa|^{\alpha}\, t^{\beta}}{\Gamma(1+\beta)}+
\frac{|\kappa|^{2\alpha}\, t^{2\beta}}{\Gamma(1+2\beta)}-\cdots \,, $$
 from which  for $t>0$ we get 
 $$ 
 \sigma^2(t)= 
 \cases{
 {\ds \,2 \,t}\,,  & $\; \alpha=2\,,\; \beta=1 \q \hbox{(normal diffusion)}\,,$ \cr
 {\ds \frac{2 t^{\beta}} {\Gamma(1+\beta)}}\,,  & $\; \alpha=2\,,\; 0<\beta<1 \q \hbox{(sub-diffusion)}\,,$ \cr
  \infty\,,  & $ \; 0<\alpha<2\,,\; 0<\beta\le 1  \q \hbox{(super-diffusion)}\,.$
  }
\eqno(2.10)$$
Admitting that these calculations are formal we leave the task 
of strict justification to pure mathematicians.
\section{The continuous time  random walk} 
In the Sixties and Seventies of the past century
 Montroll and  Weiss and Scher 
 (just to cite  these pioneers)
 published a  series
of papers for modelling
rather general diffusion processes by random walks on lattices, see \eg
 \cite{MontrollScher 73,MontrollWeiss 65}, and
the book by Weiss \cite{Weiss BOOK94} with  references therein.
Initiated by their activity the concept of continuous time random walk 
became popular in physics. 
CTRWs are good phenomenological models for several types of diffusion,
in particular the {\it microscopic} aspects of a particle jumping from point to point
admitting waiting times between jumps. 
Allowing all space instead of only a lattice, a CTRW can be viewed as a 
compound renewal process or a renewal process with reward (see \cite{Cox RENEWAL67}), 
or a random walk subordinated to a renewal process. 
Let us recall  the basic notions of the CTRW theory. 
\vsp
A CTRW  
is generated by a sequence
of  independent identically  distributed ($iid$)
 positive  random  waiting times $T_1, T_2, T_3, \dots ,$
each having the same probability density function
 $\phi(t)\,,$   $\, t>0\,, $ and
a sequence of $iid$ random jumps $X_1, X_2, X_3, \dots, $
in $\RR\,,$ each having the same probability density
$w(x)\,,$ $\, x\in \RR\,.$
Setting
$t_0=0\,,$ $\, t_n = T_1+T_2 + \dots T_n$ for $n \in \NN\,,$
the wandering particle 
makes a jump
of length $X_n$ in instant $t_n$,
so that its position is $x_0=0$ for
$0\le t <T_1= t_1\,,$       and
$x_n =    X_1 + X_2 + \dots X_n\,,$
for $ t_n \le t < t_{n+1}\,. $
We require the distribution of
the waiting times
and  that of the jumps to be independent
of each other.
 We allow the probability densities  $\phi$  and $w$    to be
    generalized functions in the sense of Gel'fand  and Shilov
    \cite{GelfandShilov 64}.
\vsp
Natural probabilistic arguments lead us  to the
{\it integral  equation} for the probability density   $p(x,t)$
(a density with respect to the variable $x$)
of the particle being in point $x$ at instant $t\,, $
see  \eg
\cite{
Gorenflo KONSTANZ01,Mainardi Bonn00,Scalas PhysA05,SGM 00,Scalas PRE04},
$$
   p(x,t) =  \delta (x)\, \Psi(t)\, +
  \int_0^t  \!\!  \phi(t-t') \, \lt[
 \int_{-\infty}^{+\infty}\!\!  w(x-x')\, p(x',t')\, dx'\rt]\,dt'\,,
\eqno(3.1)  $$
 in which  the {\it survival function}
 $\Psi(t) = \int_t^\infty \phi(t') \, dt'$
denotes the probability that at instant $t$ the particle
 is still sitting in its starting position
$x=0\,. $
Clearly, (3.1) satisfies the initial condition
$p(x,0) = \delta (x)$, and  
$p$ in place of $u$ has the properties (a) and (b) of Section 1. 
In the Fourier-Laplace domain
Eq. (3.1)  reads as
  $$ \widehat{\widetilde p}(\kappa ,s)
 =   \widetilde {\Psi}(s) +  \widehat w(\kappa )\, \widetilde \phi(s)
   \widehat{\widetilde p}(\kappa ,s)\,, \eqno(3.2)  $$
and using $  \widetilde {\Psi}(s)  = {(1-\widetilde\phi(s)) /s}\,,$
explicitly
 $$
   \widehat{\widetilde p}(\kappa ,s)   =  {1-\widetilde\phi(s)  \over s}
 {1 \over 1- {\widehat w}(\kappa )\,{\widetilde \phi}(s)}\,.
 \eqno(3.3) $$
This  representation   is known 
as the the {\it Montroll-Weiss equation},
so named after the authors of \cite{MontrollWeiss 65},
who derived it in 1965 as the basic equation for the CTRW.
By inverting the transforms one can
    find the evolution  of the sojourn density $p(x,t)$ for
    time $t$    running from zero to $\infty$.
In fact, recalling that $|\widehat w(\kappa)| < 1$ and
$|\widetilde\phi (s)| < 1$,
if $\kappa \not= 0$ and
$s \not= 0$, Eq. (3.3) becomes
$$
\widetilde{\widehat p}(\kappa, s) =
\widetilde \Psi(s)\, \sum_{n=0}^{\infty}
[\widetilde \phi (s) \, \widehat w(\kappa)]^n =
 \sum_{n=0}^{\infty} \widetilde v_n(s)\, \widehat w_n (\kappa)
 \,,
\eqno(3.4) $$
and 
we promptly obtain the
{\it series representation of the continuous time random walk},
see \eg
 \cite{Cox RENEWAL67}  (Ch. 8, Eq. (4))
or \cite{Weiss BOOK94}  (Eq.(2.101)),
$$
p(x,t) = \sum_{n=0}^{\infty} v_n(t)\, w_n (x)
=    \Psi(t)  \, \delta (x) + \sum_{n=1}^{\infty} v_n(t)\, w_n (x)\,,
\eqno(3.5)$$
where the functions  $v_n(t)$
and $w_n(x)$ are obtained by iterated convolutions
in time  and in space,
$ v_n(t) = (\Psi * \phi^{*n})(t)$,
and $\, w_n(x) = (w ^{*n})(x)$, respectively.
In particular,
$ v_0(t) = (\Psi*\delta) (t)= \Psi(t)$, $\, v_1(t) = (\Psi * \phi)(t)$,
$\, w_0(x) = \delta(x), \; w_1(x) = w(x).$
In the R.H.S  of Eq (3.5) we  have isolated
the first singular term  related to the initial condition
$p(x,0) =  \Psi(0) \, \delta (x) =\delta (x)$.
The representation (3.5)
can be found without detour over (3.1) by direct probabilistic reasoning
and transparently exhibits the CTRW  as   a  {\it subordination} of a random
    walk to a {\it renewal process}:
it can be used as starting point to derive the
Montroll-Weiss equation, as it was originally recognized by Montroll and
Weiss \cite{MontrollWeiss 65}.
\vsp
With the  special choice 
$\phi (t) = m\e^ {-mt}$, $m>0$,
  equation (3.1) describes
the {\it compound Poisson process}.  
It reduces after some
manipulations (best  carried out in the Fourier-Laplace domain)
  to the {\it Kolmogorov-Feller equation}:
$$
  \frac{\d }{\d  t}\,p(x,t)= -m\, p(x,t)+ m\,\int_{-\infty}^{+\infty}
w(x-x')\, p(x',t)\, dx' \, .
\eqno(3.6)
$$
Then from (3.5) we  obtain the    series representation 
$$p(x,t) = \sum_{k=0}^{\infty} \frac{ (mt)^k}{k!} \,
{\e}^{- mt} \,w_k (x)\,.
\eqno(3.7)
$$
Note that only in this case the corresponding stochastic process
 is {\it Markovian}.
\section{Relevance of power laws and well-scaled passage to the diffusion limit}
In this section we work out the effect the power laws for the distributions of jumps and waiting times 
have on the limiting properties of continuous time random walks. 
By power law we usually  mean  a law of decaying at infinity like a power
(with negative exponent) of the relevant  independent variable.  
There is a vast amount of literature on such laws. Let us recommend just two items, 
namely  \cite{Newman 05} and \cite{Schroeder 91}. 
What we are going to show now is that appropriate power laws 
for jumps and waiting times 
are microscopic models for fractional diffusion processes.
\vsp
     To be sufficiently general we introduce the cumulative functions for  waiting times 
and jumps. With our densities $\phi(t)$  and $w(x)$   we define
$$
\Phi(t) \!= \! \int\limits_{0}^{t}\phi(t')\, dt', \; 0 \le t < \infty; \; 
   W(x)\!= \!\int\limits_{-\infty}^{x}w(x')\,dx',
  \; -\infty <x<\infty.
 \eqno(4.1)$$
These functions are non-decreasing, 
$\Phi(0)\!= \! W(-\infty) \!= \! 0,$  $\Phi(\infty) \! = \! W(\infty) \! = \!1.$
\\    
 As they may have 
points of discontinuity we agree on the provision that equations in which they occur are meant to 
hold at points of continuity.
	 The notion of {\it power law}  concerns the behaviour of $\Phi(t)$    and $W(x)$   for large $t$  and large 
 $|x|$, roughly the mode of decrease  of the {\it tails}
 $1-\Phi(t)= \Psi(t)$ near $t=\infty$ and  $W(x)$      near $x=-\infty$,
 $1-W(x)$   near $x=+\infty$, like a power of  $t$  or $|x|$  with a negative exponent. 
  \vsp  
	 Not wanting to overload our presentation we assume  {\it spatial symmetry} 
	 ($\theta=0$) with respect to  $x=0$  and  avoid decoration of asymptotic  behaviours with slowly varying 
functions. For such neglected generalizations we refer to 
\cite{GAR Vietnam04,GorMai CARRY04}.
	The parameters $\alpha$   and $\beta$ 
	 of equation (2.1)  play  essential roles yielding 
	  power laws in the strict sense if 
$0<\alpha < 2$,  $0<\beta < 1$,
 but still formal  analogies in the extreme cases $\alpha = 2$, $\beta = 1$. 
  \vsp 
	 The question of interest  is the {\it long-time}, {\it wide-space} behaviour of a CTRW under power law 
assumptions for waiting times and jumps, i.e. the appearance of such CTRW when observed after 
long time and from far away. To bring the distant future and the far-away space into near sight   
we multiply  time intervals by a small positive factor $\tau$, space intervals by a small positive factor $h$, 
so  making  large intervals {\it numerically} of moderate size, intervals of moderate size {\it numerically} small. 
Essentially this means changing the units of time and space from 1 to $1/\tau$  and
$1/h$, respectively. 
We then obtain the asymptotic behaviour  by sending $\tau$   and  $h$  to zero, in a specially combined 
way, namely under the requirement of a  {\it scaling relation}, honouring which we call the whole 
procedure  "{\it well-scaled passage to the diffusion limit}". 
In this limit we will obtain a process 
obeying the space-time fractional diffusion equation (2.1) with $\theta=0$ , a fact that justifies the 
CTRW approach.  Conversely, we can consider a CTRW as a model  
"in the small" of a {\it space-time fractional  diffusion} process.
For generalization to skewed processes ($\theta \ne 0$) see \cite{GorMai CARRY04}. 
\vsp
    As we carry out the essential work in the Fourier-Laplace domain we use the fact that the 
behaviour of functions  $f(t)$, $g(x)$ in the infinite is mirrored in that of their transforms  
$\widetilde f(s)$, $\widehat g(\kappa)$   for $s\to 0^+$, $\kappa \to 0$. 
The lemmata we need are provided by the Tauber-Karamata 
theory for which we refer to 
\cite{Bingham-et-al RV87} and \cite{Feller 71}, 
they can also be distilled from the Gnedenko 
theorem on the domains of attraction of stable probability laws, see 
\cite{GnedenkoKolmogorov 54}. See \cite{GAR Vietnam04} and \cite{GorMai CARRY04}
for more general versions.
What we need is contained in the following two lemmata
which are simplified versions of more general facts. 
For the reader's convenience we give proofs in the Appendix. 
\newpage
\noindent
{\bf MASTER LEMMA 1:}
\\
\it
 Assume  $W(x)$ increasing: $W(-\infty)\!=\!0, W(+\infty)\!=\!1,$
 \\
 symmetry:
 $\int\limits_{(-\infty,-x)} \!\! dW(x')=\int\limits_{(x,+\infty)} \!\! dW(x')$
   for $x \ge 0$, and
  either (a) or (b):
\\
(a)  $\sigma^{2}:=\int\limits_{(-\infty, +\infty)} x^{2} dW(x) <\,\infty$,
labelled as  $\alpha=2$ ,
\\
(b) $\int\limits_{(x,\infty)}d W(x') \sim b \alpha^{-1} x^{-\alpha}$ 
  for $x \to +\infty$,  $\alpha \in (0,2)$  and   $b>0$. 
\\
  Then we have the asymptotics  $1-\widehat{w}(\kappa ) \sim
\mu|\kappa|^{\alpha}$ for $\kappa \to 0$  
  with
$\mu={\sigma^{2}}/{2}$  {in case (a)} and 
$\mu= {b \pi}/{[\Gamma(\alpha +1) \,\sin (\alpha \pi /2)]}$  
{in case (b)}.
%
\vsp
\rm
\textbf{MASTER LEMMA 2:}
\\
\it
Assume $\Phi(t)$ increasing: $\Phi(0)=0$, $\Phi(+\infty)=1$,  and  either (A) or (B): 
\\
(A) $\rho:=\int\limits_{(0, +\infty)} t \,d \Phi(t)< \infty $, labelled as
$\beta=1$, 
\\
(B) $ \Psi(t)=\int\limits_{(t,\infty)} d \Phi(t)\sim c \beta^{-1} t^{-\beta}$ for $t\to \infty$, $\beta \in (0,1)$ 
and  $c>0$.
\\
 Then
 we have the asymptotics $ 1- \widetilde{\phi}(s) \sim \lambda s^{\beta}$ for $0<s \to 0$
 with
$\lambda=\rho$ {in case (A)} and
$\lambda = c \pi/ [\Gamma(\beta+1)\, \sin(\beta \pi)]$
{in case (B)}.
\vsp
 \rm
     Assuming the conditions of these two lemmata satisfied we are ready for passing to the 
diffusion limit. We multiply the jumps $X_k$  by  a factor $h>0$, 
the waiting times $T_k$  by a factor $\tau>0$. 
 So we get a transformed random walk 
 $S_{n}(h,\tau)=\sum\limits_{k=1}^{n} h X_{k}$  with jump instants
  $t_{n}(\tau) =\sum\limits_{k=1}^{n} \tau T_{k}$  that we now
  investigate with the aim of passing to the limit 
  $h\to 0$, $\tau \to 0$ under a {\it scaling relation} between $h$ and $\tau$  yet to be established.
  As it is convenient to work in the
  Fourier-Laplace domain we note that the density $\phi_{\tau}(t) $
  of the reduced waiting times $\tau T_{k}$   and the density
  $w_{h}(x)$ of the reduced jumps $h X_{k}$  are $\phi_{\tau}(t)=\phi(t/ \tau)/ \tau $,
  $t \ge 0$  ; $w_{h}(x)=w(x/h)/h $, $-\infty <x< \infty$. The
  corresponding transforms are simply $\widetilde{\phi_{\tau}}(s)=
  \widetilde{\phi}(s \tau)$,
  $\widehat{w_{h}}(\kappa)=\widehat{w}(\kappa h)$. We are interested
  in the sojourn probability density $p_{h,\tau}(x,t)$   of the
  particle subject to the transformed random walk. In analogy to the
  Montroll-Weiss equation (3.3) we get
$$
\widehat{\widetilde{p}}_{h,\tau} (\kappa,s) =\frac{1-\widetilde{\phi}_{\tau}(s)}{s}
\frac{1}{1-\widehat{w}_{h}(\kappa) \widetilde{\phi}_{\tau}(s)} \;= 
\frac{1-\widetilde{\phi}(\tau s)}{s}
\frac{1}{1-\widehat{w}(h\kappa) \widetilde{\phi}(s\tau)} \, .
\eqno(4.2)
$$
Considering now $s$    and $\kappa $ \emph{fixed} and $\ne 0$ we find for $h\to 0$,
$\tau \to 0$  from the Master Lemmata (replacing there $\kappa$   by
$\kappa h$,  $s$  by  $s \tau$ )  by a trivial calculation, omitting
asymptotically negligible terms, the asymptotics (4.3) with (4.4).
$$
\widehat{\widetilde{p}}_{h,\tau}(\kappa,s) \sim \frac{\lambda
  \tau^{\beta} s^{\beta-1} }
  {\mu (h|\kappa|)^{\alpha}+ \lambda (\tau s)^{\beta} } = 
  \frac{s^{\beta -1}}{r(h, \tau)|\kappa|^{\alpha} + s^{\beta} } \, ,
\eqno(4.3)$$
$$
 r(h,\tau)= \frac{\mu h^{\alpha}}{\lambda\tau ^{\beta} } \, .
\eqno(4.4)$$
So we see that, for every fixed  real $\kappa \ne 0$  and  positive $s$,
$$
\widehat{\widetilde{p}}_{h,\tau}(\kappa,s) \,\to \,\frac{s^{\beta   -1}}{|\kappa|^{\alpha}+ s^{\beta}
}=\widehat{\widetilde{u}}(\kappa,s) \, ,  
\eqno(4.5)$$
as $h$ and $\tau$ tend to zero under the \textit{scaling relation}
$r(h,\tau)\equiv 1$.
Comparing with (2.7) we recognize here
$\widehat{\widetilde{u}}(\kappa,s) $  as the combined Fourier-Laplace
transform of the solution to the Cauchy problem
(2.1) with $\theta=0$. Invoking now the continuity theorems of probability
theory (compare \cite{Feller 71}),
bypassing some analytical subtleties, we see that the  time-parametrized  sojourn
probability density converges {\it weakly} (or {\it in law}) to the solution of
the Cauchy problem (2.1) with $\theta =0$.
This weak convergence can be taken as justification for approximate simulation 
of trajectories (particle paths) by CTRW's with power law jumps and waiting times 
so chosen that routines are available for producing the needed random numbers 
(see e.g. \cite{GAR Vietnam04}).
 \section{The Mittag-Leffler waiting time law and time-fractional processes}
We now offer another view to the transition from CTRW (under power law assumptions) 
to fractional diffusion, separating the temporal and spatial limiting procedures,
 thereby getting among other matters the time-fractional 
 CTRW discussed in the pioneering paper \cite{HilferAnton 95}. 
\vsp
Turning our attention to time-fractional processes we present in a condensed way some results from our papers
\cite{GorMai BAD-HONNEF07,Gorenflo KONSTANZ01,Mainardi Bonn00}. 
  We will see the importance of the Mittag-Leffler waiting time density  
 $\phi^{ML}(t)$ and the corresponding survival function $\Psi^{ML}(t)$  with their 
 Laplace transforms  displayed here:
 $$
 \phi^{ML}(t)= -(d/dt)E_\beta (-t^\beta)\,, \q  \widetilde \phi^{ML} (s) = {\ds \rec{1+s^\beta }}\,,
 \eqno(5.1)$$
 $$ \Psi^{ML}(t) = E_\beta(-t^\beta)\,, \q  \widetilde \Psi^{ML} (s) = {\ds \frac{s^{\beta-1}}{1+s^\beta }}\,,
 \eqno(5.2)$$
 with $E_\beta$  defined in (2.9).
 Throughout we assume  $0<\beta \le 1$. We recall that for $\beta=1$
 we recover $\phi^{ML}(t)= \exp (-t)$, $\Psi^{ML}(t) = \exp (-t)$.
\vsp
The importance of these functions cannot be overestimated, they also play an essential role in the theory of 
fractional relaxation (see e.g. \cite{GorMai CISM97}). 
The relevance of the Mittag-Leffler waiting time law for time-fractional CTRW has been put in bright 
light by Hilfer and Anton in 1995 \cite{HilferAnton 95}. 
Fulger, Scalas and Germano \cite{Scalas-et-al PRE08} pay special attention to its use 
as waiting time law in CTRW simulation .
In the Sixties of the past century it 
has been found by Gnedenko and Kovalenko \cite{GnedenkoKovalenko_QUEUEING68} as a limiting law in thinning
(rarefaction) of a renewal process under the power law assumptions of our Master Lemma 2, 
but they only gave its Laplace transform without identifying it as a function of Mittag-Leffler type. 
We will show, without invoking the concept of thinning, that it is a universal limiting law for long-time behaviour 
of a renewal process under a power law regime. 
\vsp
For a general CTRW we will show how via the concept of the {\it memory function} 
we can separate the passages of the scaling factors $\tau$ and $h$ (of the preceding  Section 4) to zero, 
thus avoiding the simultaneous use of the continuity theorems for the transforms of Laplace and Fourier.
\vsp
 \subsection{Manipulations: rescaling and respeeding}  
  To introduce the memory function $H(t)$ and explain its meaning we recall the  CTRW theory of Section 3, 
  in particular equations (3.1) and (3.3). 
  We need this general theory already for embedding into it the renewal theory.  
  In fact: we can view a pure renewal process as a special type of CTRW, namely one in which all jumps are of 
  fixed size 1 and the position $x$ of the wandering particle in space plays the role of the counting number 
  $n$ of the renewal process. 
  However, we prefer to continue working in the general CTRW context.
  Introducing formally in the Laplace domain the auxiliary function 
$$ \widetilde{H} (s) =
\frac{1- \widetilde{\phi}(s) }{ s\, \widetilde{\phi}(s)}
= \frac{\widetilde{\Psi}(s) }{\widetilde{\phi}(s)}
\,, \q \hbox{hence}\q
   \widetilde{\phi}(s) =  \frac{1}{1+s \widetilde{H} (s)} \,,
 \eqno(5.3)$$   
 and assuming that its Laplace inverse $H(t)$ exists, we get,
 following \cite{Mainardi Bonn00}, 
in the Fourier-Laplace domain the  equation 
$$  \widetilde{H} (s) \, \left[
  s\widehat{\widetilde p}(\kappa ,s)-1\right] =
  \left[ \widehat w(\kappa )-1\right]\,
   \widehat{\widetilde p}(\kappa ,s)
\,, \eqno(5.4)$$
and in the space-time  domain  the generalized Kolmogorov-Feller equation
$$ \int_0^t   H(t-t')\,
 \frac{\d}{\d t'} p(x,t')\, dt'  =
      -p(x,t) + \int_{-\infty}^{+\infty} w (x-x')\,p(x',t)\,dx',
\eqno(5.5)$$
with $p(x,0) =\delta (x)$. 
With the special choice of the {\it power law memory function}
$$  H^{ML}(t)= 
\cases{
{\ds \frac{t^{-\beta}} {\Gamma(1-\beta)}} \,, & $if \q 0<\beta<1\,,$\cr
{\ds \delta(t)}\,, & $if  \q \beta=1 \,,$
}
\eqno(5.6)
$$
whose Laplace transform is
$$\widetilde H^{ML}(s)= s^{\beta -1}\,, \q 0<\beta \le 1\,,\eqno(5.7)
$$
we have the {\it Mittag-Leffler waiting time law} given by
Eqs (5.1) and (5.2).                        
In the extremal case $\beta = 1$ this reduces to the {\it exponential waiting time law}
$\phi (t) = \exp(-t)$, $\Psi(t) = \exp( -t)$,
and we obtain in the Fourier- Laplace domain
$$ 
  s\widehat{\widetilde p}(\kappa ,s)-1 =
  \lt[ \widehat w(\kappa )-1\right]\,
   \widehat{\widetilde p}(\kappa ,s)
\,, \eqno(5.8)$$
 in the space-time domain 
the {\it classical Kolmogorov-Feller equation}
$$ \frac{\d}{ \d t} p(x,t) = - p(x,t)  +
   \int_{-\infty}^{+\infty}   w  (x-x')\,p(x',t)\,dx'\,,
     \quad p(x,0) = \delta(x)\,.
\eqno(5.9) $$
For    $0 < \beta < 1$  we have  the {\it time-fractional
Kolmogorov-Feller equation}  
$$   \, _tD_*^\beta \,  p(x,t) =
     -  p(x,t) +   \int_{-\infty}^{+\infty} w(x-x')\,
   p(x',t) \, dx'\,, \;p(x,0^+) = \delta(x)\,.   \eqno(5.10)$$
 Let us now consider two types of manipulations on  the CTRW:
 \\
{\bf A: rescaling, B: respeeding}.
\vsp
{\bf (A)}  means, as in Section 4,  change of the unit of time. With the {\it positive scaling factor}
$ \tau$ we replace the waiting time $T$ by $\tau T$, again intending  $\tau << 1$. 
In a moderate span of time we will so have a large number of jump events. Again we get the 
rescaled waiting time density and its corresponding Laplace transform as
$\phi _\tau (t) = \phi (t/\tau )/\tau$, 
$\widetilde\phi _\tau (s) = \widetilde\phi (\tau s)$.
By decorating also the density $p$ with an index $\tau$ we obtain this rescaled CTRW integral equation 
in the Fourier-Laplace domain as
$$ \widetilde{H}_\tau  (s) \, \lt[
  s\widehat{\widetilde {p}}_\tau(\kappa ,s)-1\rt] =
  \lt[ \widehat w(\kappa )-1\rt]\,
   \widehat{\widetilde {p}}_\tau (\kappa ,s)\,,
\eqno(5.11)$$
where  
$$ \widetilde{H}_\tau (s) =
\frac{1- \widetilde{\phi}(\tau s) }{ s\, \widetilde{\phi}(\tau s)}
\,, \q \hbox{hence}\q
   \widetilde{\phi}(\tau s) =  \frac{1}{1+s \widetilde{H}_\tau (s)} 
\,.
\eqno(5.12)$$      
\vsp
{\bf Remark:} Note that in this Section 5 the position  of the indices 
at the density $p$
and their meaning  
are convenient but different from those in the other Sections. 
\vsp
{\bf (B)}  means multiplying the quantity representing ${\ds \frac{\d}{\dt} p (x,t)}$
by a factor $1/a$, where $a>0$ is the {\it respeeding factor}:
$a>1$ means {\it acceleration}, $0<a<1$  means {\it deceleration}.
In the Fourier-Laplace representation this means multiplying the RHS of Eq. (5.11) by the factor $a$  
since the expression $\lt[s\widehat{\widetilde {p}}(\kappa ,s)-1\rt]$,
in view of $p(x,0)= \delta(x)$ and $\widehat \delta(\kappa)=1$, 
    corresponds to
$ {\ds \frac{\d}{\dt} p (x,t)}$.
\vsp
We now consider the procedures {\bf (A)} and {\bf (B)} 
in their combination so that  in the transformed domain 
 the  rescaled and respeeded process has the form
$$  \widetilde{H}_\tau  (s) \, \lt[
  s\widehat{\widetilde {p}}_{\tau,a}(\kappa ,s)-1\rt] =
 a\, \lt[ \widehat w(\kappa )-1\rt]\,
   \widehat{\widetilde {p}}_{\tau,a} (\kappa ,s)\,.
   \eqno(5.13)$$
Clearly, the two manipulations  can be discussed separately:
the choice $\{\tau >0,\, a=1\}$ means {\it pure rescaling},
the choice $\{\tau =1,\, a >0\}$ means {\it pure respeeding}
of the original process.
In the special case  $\tau =1$   we only respeed the
system; if  $0 <\tau \ll 1$
  we can counteract the compression effected by rescaling to again obtain a
moderate number of events in a moderate span of time by respeeding
(decelerating) with  $0<a \ll 1$.
 These vague notions will become clear as soon as we consider power law
waiting times.
Defining now
$$ \widetilde{H}_{\tau,a}  (s) := \frac{\widetilde{H}_{\tau}  (s)}{a}=
\frac{1- \widetilde{\phi}(\tau s) }{as\, \widetilde{\phi}(\tau s)}\,.
\eqno(5.14)$$
we finally get, in analogy to (5.11), the equation
$$ \widetilde{H}_{\tau,a}  (s) \, \lt[
  s\widehat{\widetilde {p}}_{\tau,a}(\kappa ,s)-1\rt] =
  \lt[ \widehat w(\kappa )-1\rt]\,
   \widehat{\widetilde {p}}_{\tau,a}(\kappa,s)\,.
   \eqno(5.15)
$$
What is the combined effect of rescaling and respeeding on the waiting
time density?
In analogy to (5.3) and taking account of (5.14)  we find
$$   \widetilde{\phi}_{\tau ,a}(s) =
        \frac{1}{1+s \widetilde{H}_{\tau ,a} (s)}=
        \frac{1}
{1+s {\ds\frac{1-\widetilde{\phi}(\tau s)}{as \,\widetilde{\phi}(\tau s)}}}
\,,\eqno(5.16)$$
and so, for the deformation of the waiting time density, the
{\it essential formula}
$$   \widetilde{\phi}_{\tau ,a}(s) =
     \frac{a\,\widetilde{\phi}(\tau s) }
     {1- (1-a)\widetilde{\phi}(\tau s)}
\,. \eqno(5.17) $$      
\subsection{ Asymptotic universality of the Mittag-Leffler waiting time law  under power law regime}
   We now  recall the MASTER LEMMA 2 of Section 4 and assume the conditions stipulated there. 
   By using the statements of this lemma, taking
$$a = \lambda \tau^\beta\,,\eqno (5.18)$$
fixing $s$  as required by the continuity theorem of probability for Laplace transforms, 
 the asymptotics $\widetilde \phi(s) \sim 1 - \lambda s^\beta $ 
 for $0< \tau \rightarrow 0$ implies
$$    \widetilde{\phi}_{\tau ,\lambda \tau ^\beta} (s) 
=  \frac{\lambda \tau^\beta\, \left[1-\lambda \tau^\beta s^\beta + o(\tau^\beta s^\beta)\right]}
{1 - (1-\lambda \tau^\beta)\, \left[1-\lambda \tau^\beta s^\beta + o(\tau^\beta s^\beta)\right] }
  \to
\frac{1}{1 + s^\beta }\,, \eqno(5.19)$$     
corresponding to the 
density $\phi^{ML}(t)$.  
This formula expresses
{\bf the asymptotic universality of
the Mittag-Leffler waiting time law} that includes the exponential law for $\beta=1$.
 It says that our general power law waiting
time density is  gradually deformed into the Mittag-Leffler waiting time
density.   
\vsp
{\bf Remark:} 
Let us stress here the distinguished character of the Mittag-Leffler
waiting time density $\phi^{ML}(t)$ 
defined in (5.1).
It is easy to prove the identity
$$    \widetilde{\phi}^{ML}_{\tau,a} (s)
= \widetilde{\phi}^{ML} (\tau s/a^{1/\beta })
\q \hbox{for all} \q \tau >0, \q a>0\,, \eqno(5.20)$$
that  states the {\it self-similarity} of  the combined operation
{\it rescaling-respeeding}
for the Mittag-Leffler waiting time density. In fact, (5.20) implies
$ {\phi}^{ML}_{\tau,a} (t) =
   {\phi}^{ML}(t/c) /c$ with
$ c= \tau /a^{1/\beta }\,,$
 which means replacing the random waiting time $T^{ML}$
by   $ c\, T^{ML}$.
As a consequence, choosing  $a= \tau ^\beta $, we have
$$    \widetilde{\phi}^{ML}_{\tau,\tau^\beta} (s)
 = \widetilde{\phi}^{ML} ( s)
\q \hbox{for all} \q \tau >0\,. \eqno(5.21)$$
Hence {\it the Mittag-Leffler waiting time density  is invariant against
    combined rescaling  with $\tau$ and respeeding with  
    $a=\tau^\beta$}.
	Observing  (5.19) we can say that $\phi^{ML} (t)$ is a $\tau \to 0$
	attractor for any power law waiting time (compare Master Lemma 2) with
	$$  
 \Psi(t) \sim \frac{c}{\beta} \, t^{-\beta } \,,
 \q 0<\beta<1\,,\q c>0\,, \eqno(5.22)$$
	 under combined rescaling with $\tau$ and respeeding with
	$a= \lambda \tau^\beta$.
This attraction property of the Mittag-Leffler waiting time 
 distribution with respect to power law  waiting  times (with $0<\beta <  1$) 
 is a kind of analogy to the attraction of sums of power law jump 
 distributions by stable distributions.
\subsection{Diffusion limit in space}
   We can obtain from (5.5) the fractional Kolmogorov-Feller
equation  (5.10)  
for time fractional CTRW by  direct insertion of the Mittag-Leffler memory function into 
the equation or, as the previous considerations show,
 by  manipulating it  under power law assumption for the waiting time and passing to the limit. 
We have not yet operated on the jumps. To do this now, we assume the conditions of Master Lemma 1 to hold. 
Then, by another respeeding, in fact an acceleration (that we earlier had carried out in 
\cite{Gorenflo KONSTANZ01}), 
we will arrive at diffusion processes fractional in space. 
We have
{\bf three choices}: 
\\
{\bf (A)}:  {\it diffusion limit in space only, for general waiting time},
\\
{\bf (B)}:  {\it diffusion limit in space only, for Mittag-Leffler waiting time},
\\
{\bf (C)}: {\it joint limit in time and space
(with power laws in both) with
scaling relation}.
\vsp
Note that {\bf (B)}  is just a special case of {\bf (A)} but of particular relevance.
In all three cases we rescale the jump density by a factor
$h>0$, replacing the random     jumps $X$ by  $hX$.
This means changing the unit of measurement in space from $1$ to $1/h$,
 with $0 <h \ll 1$,
  so bringing into near sight the  far-away space.
The rescaled jump density is  $w_h(x) =  w(x/h)/h$,
corresponding to   $\widehat w_h(\kappa) =\widehat w_(h\kappa)$.
\vsp
{\bf Choice (A):}
Starting from the Eq. (5.4), the Fourier-Laplace representation of the CTRW equation,
 without special assumption on the waiting time density,
  we  accelerate the spatially rescaled
process by the respeeding factor $ 1/(\mu h^\alpha )$,
arriving at the equation (using  $q_h$   as new dependent variable)
$$   \widetilde{H} (s) \, \lt[
  s\widehat{\widetilde q}_h(\kappa ,s)-1\rt] =
 \frac{ \widehat w(h \kappa )-1}{\mu h^\alpha }\,
   \widehat{\widetilde q}_h(\kappa ,s)\,. \eqno(5.23) $$
Then, {\it fixing $\kappa$} as required
by the continuity theorem of probability theory for Fourier transforms,
and {\it sending $h$ to zero} we get, 
noting that
by Master Lemma 1
 $ [\widehat w(h \kappa )-1]/(\mu h^\alpha) \to -|\kappa|^\alpha$,
and
writing  $u$   in place of  $q_0$,
$$  \widetilde{H} (s) \, \lt[
  s\widehat{\widetilde u}(\kappa ,s)-1\rt] =
-   |\kappa |^\alpha
\,\widehat{\widetilde u}(\kappa ,s)\,, \eqno(5.24)$$
still with  $\widetilde H(s)$ as in (5.3).
Here  $-|\kappa |^\alpha $
  is the symbol of the Riesz pseudo-differential operator
$\, _x D_{0}^{\,\alpha}$ (known as the Riesz fractional
derivative of order $\alpha$) obtained from the Riesz-Feller fractional
derivative for  $\theta = 0$, see (2.1) and (2.4). 
We thus arrive at the integro-pseudo-differential equation 
 $$ \!\! \int_0^t   H(t-t')\,
 \frac{\d}{\d t'} u(x,t')\, dt' \, =   \,
    _x D_{0}^{\,\alpha} \,u(x,t)
      \, ,\; 0<\alpha \le 2 \,,\; u(x,0)=\delta(x)\,.\eqno(5.25) $$
{\bf Comments}:
By this rescaling and acceleration the jumps become smaller and
smaller, their number in a given span of time larger and larger,
the waiting times between jumps smaller and smaller.
In the limit there are no waiting times anymore,
the original waiting
time density $\phi (t)$   is now only spiritual, but still determines via
$H(t)$ the memory of the process.
Eq. (5.25) offers a great variety of diffusion processes with memory depending
on the choice of the function $H(t)$.
\vsp
{\bf Choice (B)}:  Inserting in (5.25) the Mittag-Leffler memory function (5.6), 
we immediately get the {\it space-time fractional diffusion equation} (2.1) with $\theta = 0$. 
\vsp
{\bf Choice (C)}: Assuming the conditions of both Master lemmata  fulfilled, rescaling as described the 
waiting times and the jumps by factors $\tau$ and $h$, starting from (5.13), 
decelerating by a factor $\lambda \tau^\beta$ in time, then accelerating 
for space by a factor $1/(\mu h^\alpha )$,
 we obtain, 
 by fixing $s$ and $\kappa $,
 the equation
$$  \widetilde{H}_\tau  (s) \, \left[
  s\widehat{\widetilde {p}}_{\tau, a(\tau ,h)}(\kappa ,s)-1\right] =
 a(\tau ,h)\, \left[ \widehat w_h(\kappa) - 1\right]\,
   \widehat{\widetilde {p}}_{\tau, a(\tau ,h)} (\kappa ,s)\,,$$
  with $\widehat w_h(\kappa) = \widehat w(h\kappa)$, 
  $a(\tau ,h) = {\lambda \tau ^\beta }/{(\mu h^\alpha) } $ and
$$
\widetilde{H}_\tau (s)=
   \frac{1- \widetilde{\phi}(\tau s) }
{ s\, \widetilde{\phi}(\tau s)} \sim \lambda \tau ^\beta s^{\beta -1 }\,,
\q \hbox{for}\q \tau \to 0\,.$$
Observing
$$\frac{ \widehat w(h \kappa )-1}{\mu h^\alpha } \to -|\kappa|^\alpha\,,
\quad \hbox{for} \quad h \to  0\,,$$
then,  
introducing the relationship of {\it well-scaledness}
$$ a(\tau,h)= \frac{\lambda \tau ^\beta}{\mu h^\alpha } \equiv 1\ \eqno(5.26)$$
between the rescaling factors  $\tau$ and $h$,  we finally get
   the limiting equation
  $$
  s^{\beta -1} \, \lt[
  s\widehat{\widetilde u}(\kappa ,s)-1\rt] =
- |\kappa |^\alpha  \,
   \widehat{\widetilde u}(\kappa ,s)\,,\eqno(5.27)  $$
 corresponding to Eqs. (2.5) and (2.1) with $\theta=0$, 
 the symmetric space-time  fractional diffusion equation.
\vsp
{\bf Comments:}
Let us point out an advantage of splitting the passages
$\tau \to 0$ and $h \to 0$.
Whereas by the combined passage as in choice (C), if done in the well-scaled way (5.26),
the mystical concept of respeeding can be avoided, there arises
the question of correct use of the continuity theorems
of probability. There is one continuity theorem for the Laplace transform,
one for the Fourier transform, see \cite{Feller 71}.
Possible doubts whether their simultaneous use is legitimate
vanish by applying them in succession, as in our
two splitting methods.
For a more detailed discussion of mathematical aspects and the involved stochastic processes we refer to 
our recent paper \cite{GorMai BAD-HONNEF07}.
\section{Subordination in stochastic processes }
The common method for simulating particle trajectories consists in interpreting 
the concept of subordination (see \cite{Feller 71}) as one 
of transforming a stochastic process    $Y(t_*)$ 
(we call it {\it the parent process})  where $t_*$ is not the physical but an operational time into a 
process $X(t)$ in physical time $t$, by generating the operational time $t_*$ from the physical time $t$ 
via a positively oriented stochastic process $T_*(t)$ to arrive at the representation  $X(t) = Y(T_*(t))$. 
For simulation one then needs a routine for generating the process $T_*(t)$. 
For simulating trajectories in space-time fractional diffusion (2.1) this 
requires simulation of the hitting time process 
which is inverse to the stable subordinator in Feller's  parametrization
\cite{Feller 71}, 
the stable process of order $\beta$ and skewness ${-\beta}$.
\vsp   
There are routines available for simulating stable variates, see e.g. \cite{Janicki LN96,Janicki-Weron 94}.
But we do not know of easy routines for inverting a stable subordinator. 
Therefore, we have developed our method of {\it parametric subordination} which, 
by starting from the operational time $t_*$ allows construction of trajectories 
by  producing: 
FIRST an  $\alpha$-stable process $x = Y(t_*)$ with skewness $\theta$ for 
for the position $x$ in space, 
 SECOND an extreme positive-oriented  stable process $t = T(t_*)$  
of order $\beta$  with skewness  $-\beta$ that we call  the {\it leading process}.
Then we get, in the $(t, x)$-plane, the parametrized graph  $t = T(t_*),  x = Y(t_*)$ of a desired trajectory 
of the process  $X(t)$  corresponding to eq. (2.1)  as $X(t) = Y(T_*(t))$ where now   $t_* = T_*(t)$  
is the process 'inverse' 
 to  $t = T(t_*)$. For the more general situation we 
 refer to our recent paper \cite{Gorenflo-Mainardi-Vivoli CSF07}.
 This method is exact in the sense that it allows us to produce 
 {\it snapshots of a true particle path}. 
 Let us in this context also draw attention to the recent paper 
 \cite{Friedrich 07}  by Kleinhaus and Friedrich. 
\vsp
   Let us sketch  how this method directly arises from the CTRW approximation 
under appropriate power law assumptions on the waiting time and the jumps. 
Compare also with \cite{Saichev PhysA05} and \cite{SokolovKlafter EINSTEIN05}.
For mathematical 
details we refer to \cite{Gorenflo-Mainardi-Vivoli CSF07}. 
We start with the equations (3.4) and (3.5), 
and rescale waiting times and jumps 
again with positive factors $\tau$ and $h$. 
In the Fourier-Laplace domain we assume,  a bit more general than in our Master lemmata, 
the power-law conditions
$$ 1 - \widetilde {\phi} (s) \sim \lambda s^\beta \,, \q 
  \lambda >0 \q \hbox{as} \q  s \to 0^+\,, \eqno(6.1)$$
$$1-\widehat{w}(\kappa ) \sim \mu
|\kappa| ^\alpha \, i^{\ds \, \theta \sgn \kappa }\,,
\q  \mu >0, \q \hbox{as}\q  \kappa \to 0\,. \eqno(6.2)
$$
We see that we have a combination of two Markov processes happening in discrete time, one giving a jump 
(with density $w(x)$)
 in space $x$ at every instant $n$ where $n$ is the running index,  
 the other one giving a positive jump in time $t$ at every instant $n$.  
\vsp
Using the effect of the rescaling on Eq. (3.5)  and correspondingly 
decorating it with additional indices 
we get  in the Fourier-Laplace domain 
$$ \widehat{\widetilde {p}}_{h,\tau }(\kappa ,s) =
\sum_{n=0}^{\infty} \frac{1-\widetilde {\phi} (\tau s)}{s}\,
\lt( \widetilde {\phi} (\tau s)\rt)^n\,
\lt( \widehat {w} (h\kappa )\rt)^n\,.\eqno(6.3)$$
Separately we treat the powers $\lt( \widetilde \phi (\tau s)\rt)^n$
and   $\lt( \widehat w (h\kappa )\rt)^n$, so avoiding
the problematic simultaneous inversion of the diffusion limit
from the Fourier-Laplace domain into
the physical domain.
Observing    from (6.1)
$    \lt( \widetilde \phi (\tau s)\rt)^n   \sim
    \lt( 1 - \lambda (\tau s)^\beta \rt)^n $,
we relate the running index $n$ to the presumed operational time
$t_*$ by
$ n \sim {t_* }/{(\lambda \,\tau ^\beta) }$,
and for \underbar{fixed} $s$ (as required by the continuity theorem
of probability theory), by sending
$\tau \to 0$ we get
$ \lt( \widetilde \phi (\tau s)\rt)^n   \sim
\lt(1- \lambda \,\tau ^\beta s^\beta \rt)^{t_* /(\lambda \tau ^\beta )}
\to  \exp \lt(-t_*  \,s^\beta \rt)$.
Here the Laplace variable $s$ corresponds to physical time $t$, and in Laplace inversion
we must treat $t_* $ as a parameter.
Hence, in physical time $\exp (-t_* s^\beta )$ corresponds
to
$$ \bar g_\beta(t, t_*) = t_* ^{-1/\beta } \, \bar g_\beta (t_* ^{-1/\beta } t)\,,
\eqno(6.4) $$
with $\widetilde{\bar g}_\beta (s) = \exp (-s^\beta )$.
Here $\bar g_\beta(t, t_*)$  is the totally positively skewed stable
   density (with respect to the variable $t$) evolving in operational
   time $t_*$ according to the "space"- fractional equation
$$\frac{\d}{\dt_*}\, \bar g_\beta(t, t_*) =
\,_tD^\beta_{-\beta} \,\bar g_\beta(t, t_*)\,,
\quad \bar g_\beta(t,0) = \delta(t)\,,
\eqno(6.5)$$
where $t$ plays the role of the spatial variable.
Analogously, observing from (6.2)
$    \lt( \widehat w (h\kappa )\rt)^n   \sim
 \lt( 1 - \mu  (h |\kappa| )^\alpha \, i ^{\ds \theta \sgn \kappa}\rt)^n $,
 with the aim of obtaining a meaningful limit
we  set
$ n  \sim {t_* }/{(\mu h^\alpha) }$,
and find, by sending $h\to 0^+$, the relation
$$     \lt( \widehat w (h\kappa )\rt)^n
 \sim \lt(1-\mu (h |\kappa|)^\alpha\,i^{\ds \theta\sgn\kappa}\rt)
 ^{t_*/(\mu h^\alpha)}
 \to
   \exp \lt(-t_* |\kappa|^\alpha   \, i ^{\ds \theta \sgn \kappa} \rt),
     $$
the Fourier transform of a $\theta$-skewed $\alpha$-stable
density $f_{\alpha,\theta} (x,t_* )$ evolving
in operational time $t_* $.
This density is the solution of the space-fractional equation
$$ \frac{\d}{\d t_* }\,f_{\alpha,\theta} (x,t_* )
   =  \,_xD_\theta^\alpha \, f_{\alpha,\theta} (x,t_* )\,,
 \quad  f_{\alpha,\theta} (x,0)=\delta (x)\,. \eqno(6.6) $$
The two relations 
between the running index $n$ and the presumed operational time $t_* $
require the (asymptotic) {\it scaling relation}
$ \lambda \,\tau ^\beta \sim \mu \, h^\alpha$,
that for purpose of computation we simplify to
$$ \lambda \,\tau ^\beta = \mu \, h^\alpha \,. \eqno(6.7)$$
Replacing $t_* $ by $t_{*,_n} = n \lambda \tau ^\beta$, using
the asymptotic results  obtained for the powers
$ \lt( \widetilde \phi (\tau s)\rt)^n  $ and
$ \lt( \widehat w (h\kappa )\rt)^n  $, furthermore noting
$ ({1-\widetilde \phi (\tau s)})/{s} \sim
  s^{\beta -1} \,\lambda \,\tau ^\beta $,
we finally obtain from (6.3) the Riemann sum
(with increment $\lambda \tau^\beta$)
$$
 \widehat{\widetilde {p}}_{h,\tau }(\kappa ,s) \sim
s^{\beta -1} \,
\sum_{n=0}^\infty
\exp \lt[
 -n\lambda \tau ^\beta
\lt(s^\beta  +|\kappa |^\alpha i^{\ds \theta \sgn\kappa}\rt)\rt]
\,\lambda \,\tau ^\beta    \,, \eqno(6.8)
$$
and hence the integral
$$  \widehat{\widetilde {p}}_{h,\tau }(\kappa ,s) \sim
s^{\beta -1} \,
\int_{0}^\infty
\exp \lt[
 -t_* \lt(s^\beta  +|\kappa |^\alpha i^{\ds \theta \sgn\kappa}\rt)\rt]
\,dt_*\,. \eqno(6.9)$$
For the {\it limiting process} $u_\beta (x,t)$ this means
 $$
 \widehat{\widetilde {u}}_{\beta  }(\kappa ,s) =
\int_{0}^\infty s^{\beta -1} \,
\exp \lt[
 -t_* \lt(s^\beta  +|\kappa |^\alpha i^{\ds \theta \sgn\kappa}\rt)\rt]
\,dt_*\,. \eqno(6.10)$$
Observe that the RHS of this equation is just another way
of writing the RHS of equation (2.6) which is the Fourier-Laplace
solution of the  space-time fractional diffusion
equation (2.1).
By inverting the transforms we get after some manipulations
(compare \cite{M3 PRE02sub}) in physical space-time
the {\it integral formula of  subordination}
$$
u_\beta (x,t)=
  \int_0^\infty
           f_{\alpha,\theta} (x,t_*)\, g_\beta(t_*,t)\, dt_*\, \eqno(6.11)$$
  with
$$ g_\beta(t_*,t) = \frac{t}{\beta}\,
 \bar {g}_\beta \lt(t\, t_* ^{-1/\beta}\rt) \,t_*^{-1/\beta -1}\,. \eqno(6.12)  $$
Eq. (6.11) is basic for the conventional concept of subordination
where there are also two processes involved.
One is the unidirectional  motion along the
$t_*$  axis representing the operational time.
This motion happens in physical time $t$ and the $pdf$ for the
operational time having value $t_*$  is (as density in $t_*$,
evolving in physical time $t$) given by (6.12).
In fact, by substituting $y= t\, t_*^{-1/\beta }$  we find
 $$ \int_0^\infty  g_\beta (t_*,t)\, dt_* \equiv
 \int_0^\infty  \bar {g}_\beta (t,t_*)\, dt = 1\,, \quad  \forall \, t>0\,.
 \eqno(6.13) $$
The operational time $t_*$ stands in analogy to the
counting index $n$ in Eqs. (3.5) and (6.3).
The other process is the process described by Eq. (6.6), a spatial
probability density for sojourn of the particle
in point $x$ evolving in  operational time $t_*$, namely
$ \bar u_\beta(x,t_*) = f_{\alpha,\theta} (x,t_*)$.
We get the solution to the Cauchy problem (2.1), namely the $pdf$
$u(x,t)=u_\beta (x,t)$
for sojourn in point $x$, evolving in physical time $t$, by {\it averaging}
$\bar u_\beta (x, t_*)$ with the {\it weight function}
$g_\beta(t_*,t)$ over the interval $0<t_*<\infty$ according to (6.11).
\subsection{Trajectories for  space-time fractional diffusion}
In the series representation (3.5) of the CTRW the running index $n$
(the number of jumps having occurred up to physical time $t$) is a
{\it discrete operational time}, proceeding in unit steps. To this
index $n$ corresponds the physical time $t=t_n$, the sum of the first
$n$ waiting times, and in physical space the position
$x = x_n$, the sum of the first $n$ jumps, see Section 3.

We have 
{\it two discrete Markov processes} (discrete in operational time $n$),
namely a random walk in the space variable $x$, 
with jumps $X_n$, and another
random walk (only in positive direction)
of the physical time $t$, 
making a forward jump $T_n$ at every instant $n$.
Rescaling space and physical time by factors
$h$ and $\tau $, observing the {\it scaling relation} (6.7),
we introduce, by sending $ h\to 0$ and $\tau \to 0$,
the {\it continuous operational time}
$$ t_*  \sim n\,\lambda \,\tau ^\beta \sim n\, \mu \, h^\alpha \,.
\eqno(6.14)$$
Then in the diffusion limit the spatial process becomes an
$\alpha $-stable  process for the position
$x = \bar{x} =  \bar{x}(t_* )$, whereas the unilateral time
process becomes a unilateral (positively directed)
$\beta$-stable  process for the physical time
$ t = \bar{t} =   \bar{t}(t_* )$.
A trajectory of a diffusing particle in physical coordinates can
be produced by combining in the $(t,x)$ plane the two Markovian random evolutions
$$
\cases{
x = \bar{x} =  \bar{x}(t_* ) \,, \cr
t = \bar{t} =  \bar{t}(t_* ) \,, \cr
}
\eqno(6.15)
$$
obeying the stochastic differential equations (compare with \cite{Friedrich 07})
$$
\cases{
d \bar{x} =  d (\hbox{L\'evy  noise of order} \, \alpha \;
     \hbox{and skewness}\, \theta)  \,, \cr
d \bar{t} =  d (\hbox{one sided positively oriented L\'evy  noise of order} \, \beta)   \,. \cr
} \eqno(6.16) $$
This gives us in the $(t, x)$ plane
the $t_*\,$- parametrized particle trajectory that by elimination of
        $t_*$    we get  as $x=x(t)$.
We suggest to call this procedure "construction of a   particle trajectory by 
{\it parametric subordination}".
  Note that the process  $t = T(t_*)$ yielding the second
   random function in (6.16)  has the properties of
   a {\it subordinator} in the sense of Definition 21.4 in
           \cite{Sato 99}.
\vsp
\underbar{Remark}:
It is instructive to see what happens for the limiting value
      $\beta = 1$. In this case the Laplace transform of
     $\bar{g}_\beta(t,t_*) = \bar{g}_1 (t,t_*)$ is $\exp(-t_* s)$,
      implying $ \bar{g}_1 (t,t_*) = \delta(t-t_*)$, the delta density
      concentrated at $t=t_*$. So, $t=t_*$, operational
      time and physical time coincide.
\subsection{Numerical results for the symmetric case $\theta=0$}
For numerical simulation of trajectories we proceed in three steps.
\\
First, let  the operational time $t_*$ assume
 $N$ discrete equidistant values in a given interval $[0,T]$,
that is
$t_{*,n}=n T/N,\;  n=0, 1,\dots,N$.
As a working choice we take $T=1$ and $N=10^{6}$.
Then produce $N$ independent
identically distributed ($iid$) random deviates,
 $Y_1,Y_2,\dots,Y_{N}$ having a symmetric
stable probability  distribution of order $\alpha$,
see the book by Janicki \cite{Janicki LN96}
for a useful and efficient method to do that.
Now, with the  points
 $$x_0=0,\quad  x_n = \sum\limits_{k=1}^n X_k,\quad n=1,\dots,N\,,
\eqno(6.17)$$
the couples $(t_{*,n},x_n)$, plotted
in the  $(t_*,x)$ plane  (operational time, physical space)
can be considered as points of
a {\it true trajectory}
$\{x(t_*) : 0\le t_*\le T\}$
of a symmetric L\'evy motion with order  $\alpha$
corresponding
to the integer values of operational time $t_*=t_{*,n}$.
In this identification of $t_*$ with $n$ we use the fact that
our stable laws for waiting times and jumps imply
$\lambda = \mu = 1$  in the asymptotics (6.1) and (6.2) and
$\tau = h = 1$ as initial scaling factors in (6.3) and (6.7).
In order to complete the trajectory we agree to connect every
two successive  points
$(t_{*,n},x_n)$  and  $(t_{*,n+1},x_{n+1})$
by a horizontal line from $(t_{*,n},x_n)$ to
$(t_{*,n+1},x_n),$ and a vertical line from
$(t_{*,n+1},x_n)$ to $(t_{*,n+1},x_{n+1}).$
Obviously, this is not the 'true' L\'evy motion from point
$(t_{*,n},x_n)$ to point $(t_{*,n+1},x_{n+1})$, but from the
theory of CTRW we know this kind of discrete random process
to converge in the appropriate sense to  L\'evy motion.
The points $(t_{*,n}, x_n)$ are points of  a {\it true L\'evy motion},
as can be shown by observing the infinite divisibility and self-similarity of stable laws.     
\vsp
As  a second step, we  produce $N$ $\,iid$  random deviates,
  $T_1,T_2,\dots,T_{N}$ having a stable probability
distribution with order $\beta$ and skewness $ -\beta$
(extremal stable distributions).
Then, consider the points
  $$t_0=0,\quad  t_n = \sum\limits_{k=1}^n T_k,\quad n= 1,\dots,N\,,
\eqno(6.18)
$$
and plot the couples $(t_{*,n},t_n)$
in the  $(t_*,t)$ (operational time, physical time) plane.
By  connecting  points with
horizontal and vertical lines we get  snapshots of  trajectories 
$\{t(t_*) : 0\le t_*\le N \tau = 1 \}$
describing the evolution of the physical time $t$ with
increasing operational time $t_*.$
\vsp
The final (third) step consists of
plotting  points $(t(t_{*,n}),x(t_{*,n}))$ in the  $(t,x)$
plane, namely the physical time-space plane, and connecting
them as before.
So one gets a discrete   approximation of
a particle trajectory of spatially symmetric ($\theta=0$)
fractional diffusion
with  parameters $\alpha$ and $\beta$.
\vsp
Now as the successive values of $t_{*,n}$ and $x_n$ are
generated by successively adding the relevant standardized
stable random deviates,  the obtained sets of points in the
three coordinate planes:
$(t_*, t)$, $(t_*,x)$, $(t,x)$
can, in view of infinite divisibility and self-similarity of the stable
probability distributions, be considered as {\it snapshots} of
the corresponding {\it true random processes} occurring in
continuous operational time $t_*$ and physical time $t$, correspondingly.
Clearly, fine details between successive points are missing.
\begin{figure}
\begin{center}
 \includegraphics[width=.52\textwidth]{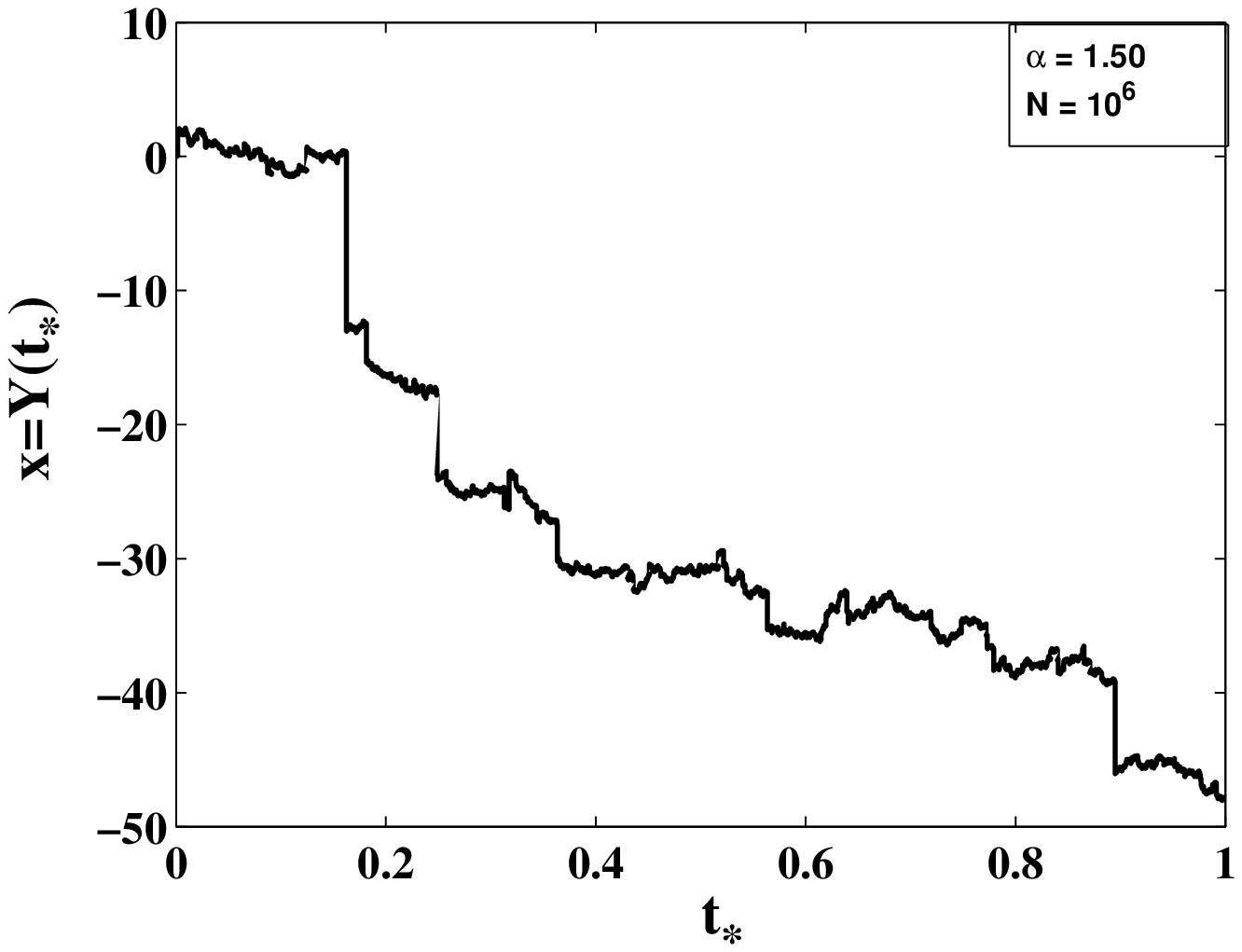}
\end{center}
 \caption{A trajectory for the parent process $x=Y(t_*)$ with
 $\{\alpha =1.5\}$.}
\vskip 0.30truecm
 \includegraphics[width=.52\textwidth]{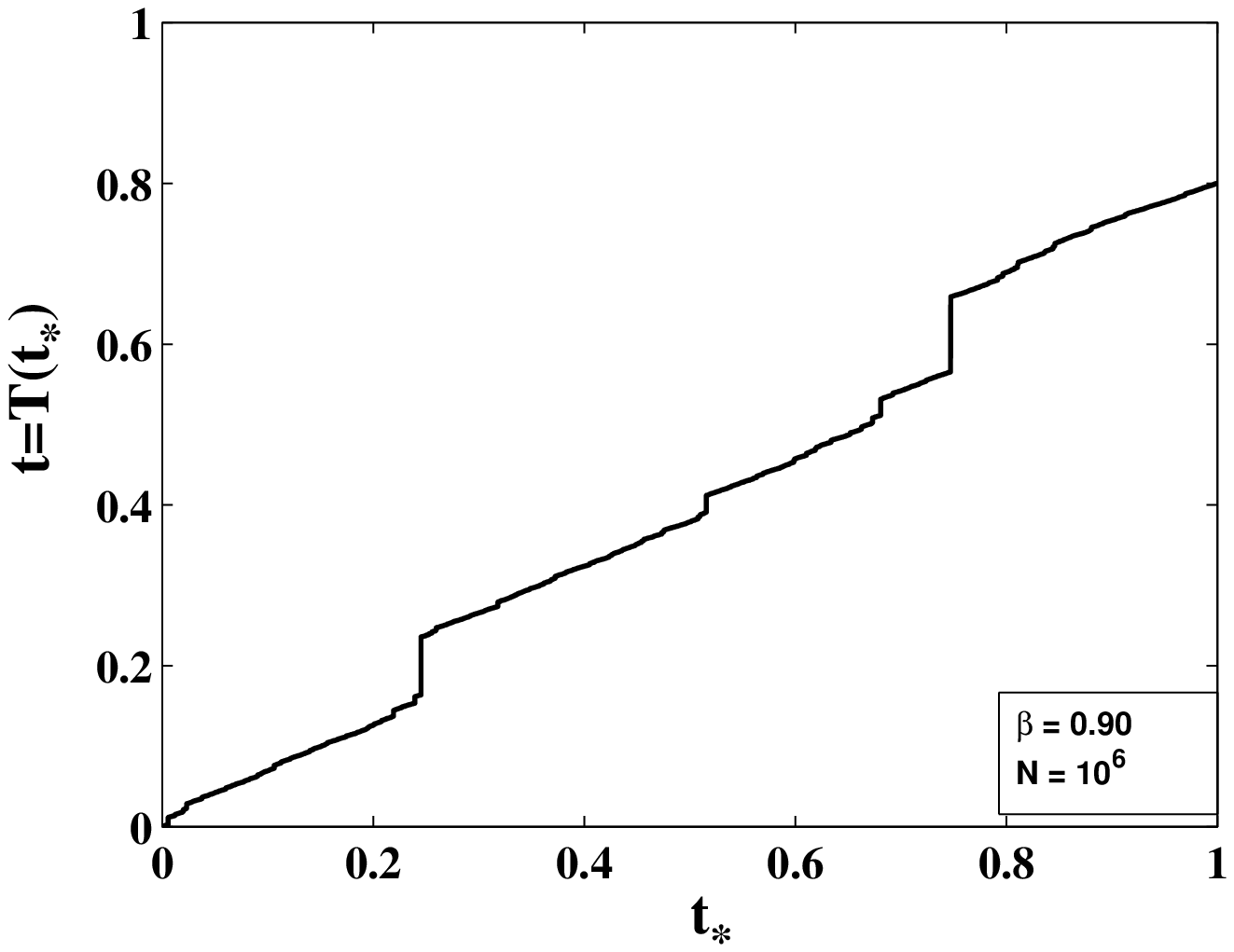}
  \includegraphics[width=.52\textwidth]{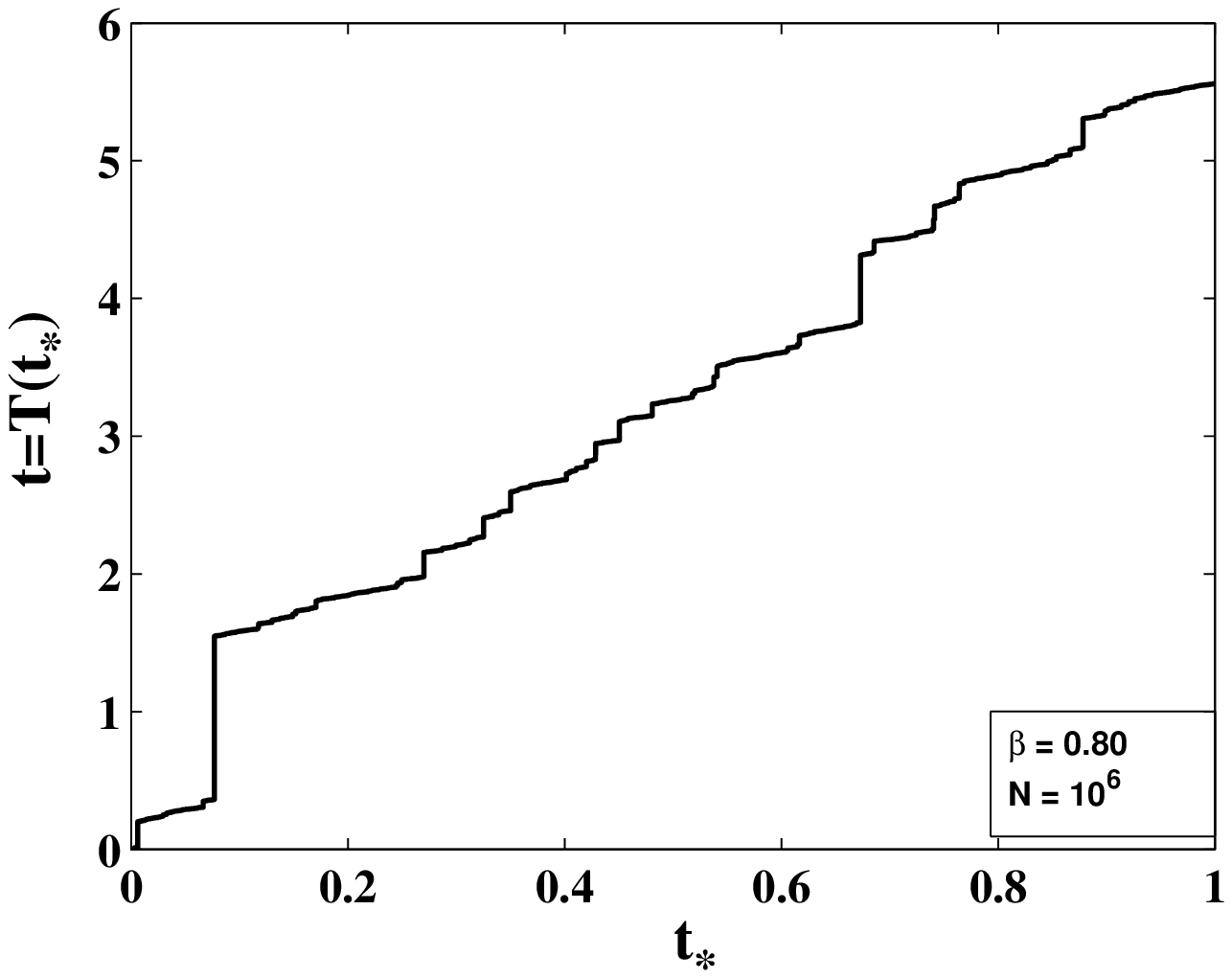}
 \caption{A trajectory for the leading process $t=T(t_*)$.}
 \centerline{LEFT: $\{\beta =0.90 \}$,
      RIGHT: $\{\beta =0.80 \}$.}
\vskip 0.30truecm
 \includegraphics[width=.52\textwidth]{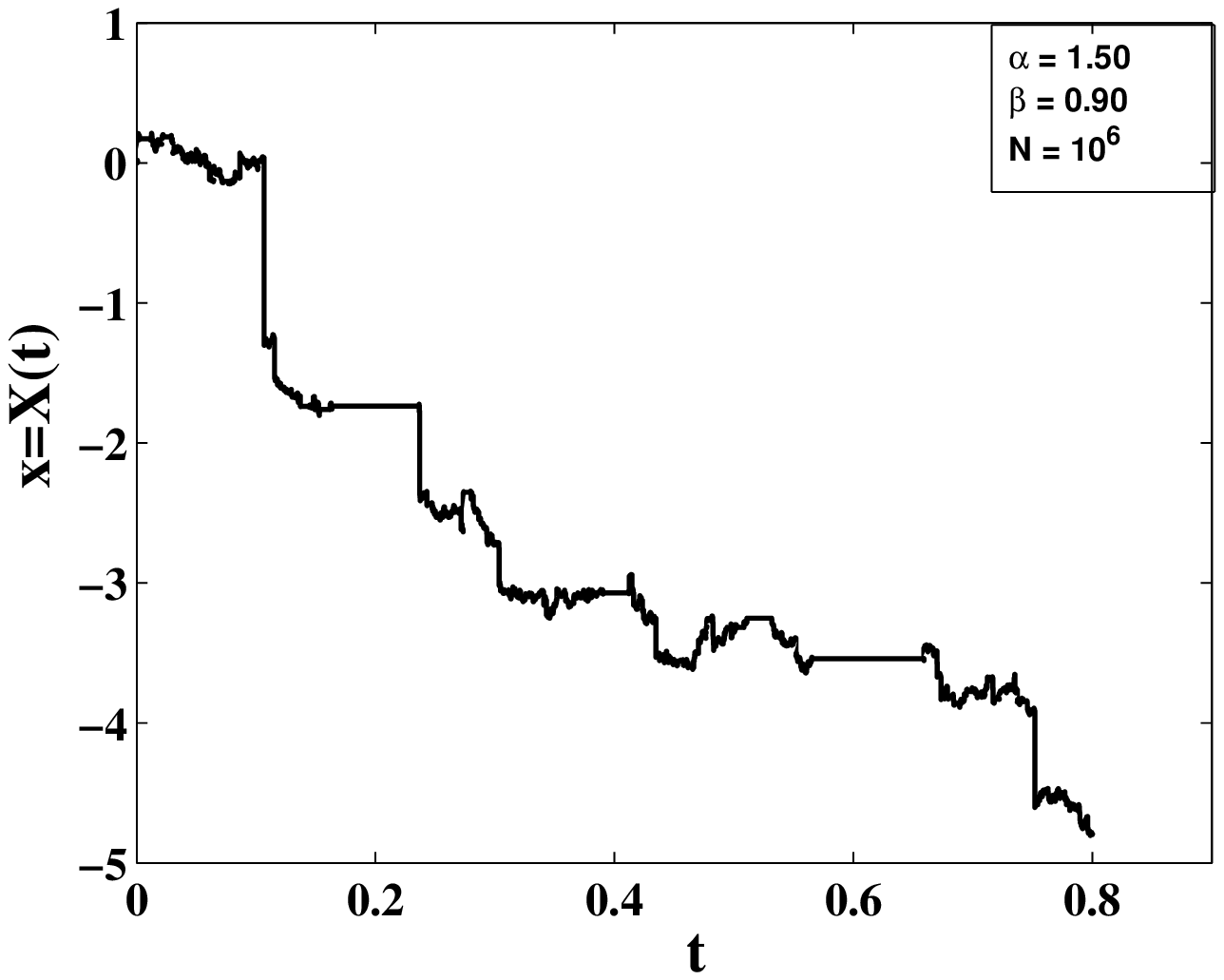}
\includegraphics[width=.52\textwidth]{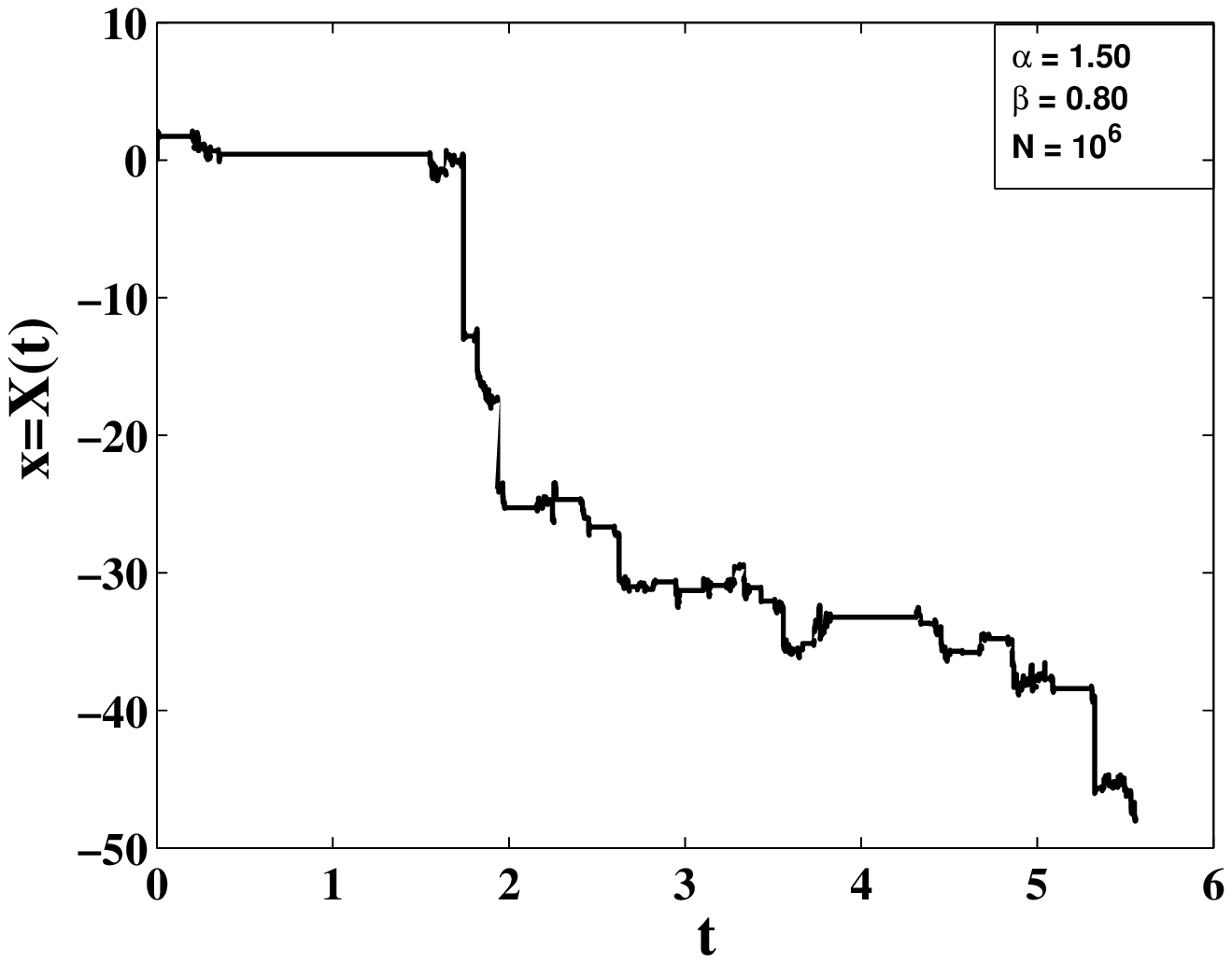}
 \caption{A trajectory for the subordinated process $x=X(t)$.}
 \centerline{LEFT: $\{\alpha =1.5,\; \beta =0.90 \}$,
      RIGHT: $\{\alpha =1.5,\; \beta =0.80 \}$.}
\end{figure}
They are hidden: 
\\
- In the $(t_*,x)$ plane in the horizontal lines from
$(t_{*,n},x_n)$ to $(t_{*,n+1},x_n)$ and the vertical
lines from $(t_{*,n+1},x_n)$ to $(t_{*,n+1},x_{n+1})$.
\\
- In the $(t_*,t)$ plane in the horizontal lines from
$(t_{*,n},t_n)$ to $(t_{*,n+1},t_n)$ and the vertical
lines from $(t_{*,n+1},t_n)$ to $(t_{*,n+1},t_{n+1})$.
\\
- In the $(t,x)$ plane in the horizontal lines from
$(t_n,x_n)$ to $(t_{n+1},x_n)$ and the vertical lines
from $(t_{n+1},x_n)$ to $(t_{n+1},x_{n+1})$.
\vsp
The well-scaled passage to
the diffusion limit here consists simply in regularly subdividing the
$\{t_*\}$ intervals  of length 1   into smaller and smaller
subintervals (all of equal length   $\tau$
and adjusting the random
increments of $t$ and $x$ according to the requirement of
self-similarity, namely taking, respectively,
the  waiting times and spatial jumps as $\tau^{1/\beta}$
multiplied by a standard extreme $\beta$-stable deviate,
$\tau^{1/\alpha}$ multiplied by a standard (in our special
case: symmetric) $\alpha$-stable deviate, respectively, as
required by the self-similarity properties of the stable
probability distributions).
Furthermore if we watch a trajectory in a large interval of
operational time $t_*,$ the points $(t_{*,n},x_n)$ and
$(t_{*,n+1},x_{n+1})$ will in the graphs appear very near to each
other  in operational time $t_*$ and aside from missing mutually
cancelling jumps up and down (extremely near to each other) we
have a good picture of the true processes.
\vsp
As interesting case studies let us present in Figs 1-3 the trajectories  
for $\{\alpha=1.5, \, \beta=0.90, \theta =0  \}$ and 
$\{\alpha=1.5, \, \beta=0.80, \theta =0  \}$, having in common the 
parent process $x=Y(t_*)$. 
\section{Conclusions}
Fractional diffusion processes as models of anomalous diffusion are gaining increasing popularity 
not only among science researchers but also among more or less pure mathematicians. 
For the latter they offer fascinating opportunities for applying pseudo-differential operators and other 
powerful analytic instruments, e.g. those of fractional calculus that in recent decades has made remarkable
 advances. 
 \vsp
 In the field of anomalous diffusion we meet challenges for the experimental sciences, 
 for mathematical modelling of real processes and their simulation, for investigation of the 
 underlying evolution processes (the macroscopic aspect) and  the fine-structure of their particle 
 trajectories (the microscopic aspect), and for numerical analysis and computational treatment 
 of less common problems. 
 In our presentation we have discussed three topics of current interest that make 
 visible the large arsenal of tools required.
\section*{Acknowledgements}
The authors appreciate the careful work of the anonymous referees. 
Their constructive critics have helped us to improve the balance and
understandability of our paper.
\section*{Appendix: Proofs of the Lemmata}
We present here proofs for Master Lemma 1 and Master Lemma 2 
of Section 4. 
We had proved analogues of these Lemmata for the probability densities 
in \cite{GorMai MAPHYSTO02}. 
Here we now give proofs for the probability distribution functions, 
thereby referring to Chapter 8 of the fundamental
treatise by Bingham et al. \cite{Bingham-et-al RV87} on Regular Variation
 through an appropriate change of notations.
 In view of this we first need to recall the notions of slowly varying functions
 and regularly varying functions. These concepts allow generalizations of the two 
 Master Lemmata stated, without proofs,  in \cite{GAR Vietnam04} and \cite{GorMai CARRY04}.
\vsp
{\bf Definitions:} We call a (measurable) positive function $a(y)$,
defined in a right neighborhood of zero, {\it slowly varying at zero} if
$a(cy)/a(y) \to 1$ with $y \to 0$ for every $c>0$.
We call a (measurable) positive function $b(y)$,
defined in a  neighborhood of infinity, {\it slowly varying at infinity}
if
$b(cy)/b(y) \to 1$ with $y \to \infty$ for every $c>0$.
An example of a slowly varying function at zero and infinity is: $|\log y|^{\gamma}$ with $\gamma \in \RR$.
Then {\it regularly
varying    functions}  are power functions multiplied by
 slowly varying functions. 
\vsp 
{\bf Proof of Master Lemma 1:}
\\
Note that because of symmetry we need only consider positive values of the 
variables $x$ and $\kappa$. 
\\
In the easy case (a) $\alpha=2$   
the well-known fact $\sigma^2 = - {\widehat w }^{\prime \prime}(0)$ 
implies 
$$1 - \widehat w(\kappa) \sim \frac{\sigma^2}{2}\, \kappa^2
\quad \hbox{as}  \quad \kappa \to 0\,.$$
\\
In case (b)    we refer to Theorem 8.1.10 in \cite{Bingham-et-al RV87}.  
It says that if for a probability distribution function $W(x)$ we set 
$$    T(x)  = W(-x) +1 - W(x)\,,  \quad
      U(\kappa) = \int_{-\infty}^{+\infty} \cos (\kappa x)\, dW(x)\,,$$
and take any function $L(x)$  slowly varying at infinity, then the relation  
 $$ T(x)  \sim L(x)\, x^{-\alpha} \quad \hbox{for} \quad x \to \infty\,,$$ 
is  equivalent to the relation 
$$ 1 -U(\kappa)  \sim  
\frac {\pi} {2\Gamma(\alpha)\, \sin(\alpha \pi/2)}
\, \kappa^\alpha \,L(1/\kappa)
\quad  \hbox{for} \quad \kappa \to 0^+ \,. $$
Taking now  $L(x)$ as the constant function 
$L(x) \equiv 2b/\alpha$ and observing that because of our 
symmetry assumption on the jump distribution function $W(x)$  we have
$W(-x)= 1-W(x)$  and hence
$ T(x) = 2 [1-W(x)] $  
in all continuity points, we arrive at $U(\kappa)= \widehat w(\kappa)$  and see that
$$ 1- W(x) = \int_x^{+\infty} dW(x') \sim b \alpha^{-1} x^{-\alpha}
\quad \hbox{for}\quad x \to \infty$$
in view of $\alpha\,\Gamma(\alpha)= \Gamma(\alpha+1)$
  implies
$$ 1- \widehat w(\kappa) \sim \frac{b \pi \kappa^\alpha}{\Gamma(\alpha+1) \, \sin (\alpha\, \pi /2)}
 \quad \hbox{for} \quad \kappa \to 0^+ \,.$$  
  We have completed the proof.
\vsp
{\bf Proof of Master Lemma 2.} 
\\
  In the easy case  (A) $\beta=1$ the statement 
  $1-\widetilde\phi(s) = \rho\,s$  is a consequence of the well-known fact
  $ \rho = -{\widetilde \phi}^{\,\prime}(0) $. 
\\
 In case (B) $0<\beta<1$ we invoke Corollary 8.1.7 of \cite{Bingham-et-al RV87}.
It says, among other things, that for a probability distribution function 
$\Phi(t)$ vanishing 
for $t<0$ the relation
$$  \Psi(t) := 1-\Phi(t) \sim \frac{1}{\Gamma(1-\beta)} \, \frac{L(t)}{t^\beta}
\quad \hbox{for} \quad t \to \infty\,,$$
where $L(t)$ is a slowly varying function at infinity, implies the  relation 
$$  1 -\widetilde \phi(s) \sim s^\beta \,L(1/s)\quad \hbox{for}\quad s\to 0^+\,.$$
 Now taking 
 $ L(t)\equiv c\,\Gamma(1-\beta)/\beta$ 
  we get
$$  1 -\widetilde \phi(s) \sim \lambda\, s^\beta \quad \hbox{for}\quad s\to 0^+\,,\quad
    \hbox{with} \quad  
	\lambda = c\,\frac{\Gamma(1-\beta)}{\beta} = 
	\frac{c \pi}{ \Gamma(\beta+1) \, \sin(\beta \pi)}\,,$$
	where  we have used  the reflection formula for the gamma function. 
The proof is complete.

\end{document}